\documentclass[times]{article}

\usepackage{amsmath,amssymb}
\usepackage[colorlinks,bookmarksopen,bookmarksnumbered,citecolor=red,urlcolor=red]{hyperref}
\usepackage{graphicx}
\usepackage{subfig}
\usepackage{stmaryrd}
\usepackage{bm}

\newcommand\BibTeX{{\rmfamily B\kern-.05em \textsc{i\kern-.025em b}\kern-.08em
T\kern-.1667em\lower.7ex\hbox{E}\kern-.125emX}}

\date{September 04, 2016}

\begin{document}

%\runningheads{M. Lawry and K. Maute}{Optimization of finite strain contact problems}

\title{Level set shape and topology optimization of finite strain bilateral contact problems}

\author{Matthew Lawry, and Kurt Maute \\ Department of Aerospace Engineering, \\ University of Colorado Boulder, 
 Boulder, CO 80309-0429, USA}

%\keywords{Topology Optimization, Level Set Methods, eXtended Finite Element Method, Stabilized Lagrange Multiplier Method, Sliding Contact, Finite Strain}

\maketitle

\section*{Abstract}
\vspace{-2pt}
This paper presents a method for the optimization of multi-component structures comprised of two and three materials considering large motion sliding contact and separation along interfaces. The structural geometry is defined by an explicit level set method, which allows for both shape and topology changes. The mechanical model assumes finite strains, a nonlinear elastic material behavior, and a quasi-static response. Identification of overlapping surface position is handled by a coupled parametric representation of contact surfaces. A stabilized Lagrange method and an active set strategy are used to model frictionless contact and separation. The mechanical model is discretized by the extended finite element method which maintains a clear definition of geometry. Face-oriented ghost penalization and dynamic relaxation are implemented to improve the stability of the physical response prediction. A nonlinear programming scheme is used to solve the optimization problem, which is regularized by introducing a perimeter penalty into the objective function. Sensitivities are determined by the adjoint method. The main characteristics of the proposed method are studied by numerical examples in two dimensions. The numerical results demonstrate improved design performance when compared to models optimized with a small strain assumption. Additionally, examples with load path dependent objectives display non-intuitive designs.

\vspace{-6pt}

\section{Introduction}
\vspace{-2pt}
Sliding contact phenomena between deformable structures play a crucial role in the functionality of many mechanical systems in commercial and industrial applications. Whether the desired functionality is to re-direct motion, provide a mechanical advantage, improve traction, or regulate stored energy, the performance of such systems is highly sensitive to interface geometry. Computational design optimization is well suited for these types of problems, as ideal design solutions can be non-intuitive. This paper provides a shape and topology optimization method for problems involving large sliding, large deformation, frictionless contact and separation in two dimensions. While interfacial adhesion and friction are ignored in this study, the proposed framework allows for the inclusion of additional contact phenomena.

\begin{figure}[t]
\begin{center}
\includegraphics[width=.5\linewidth]{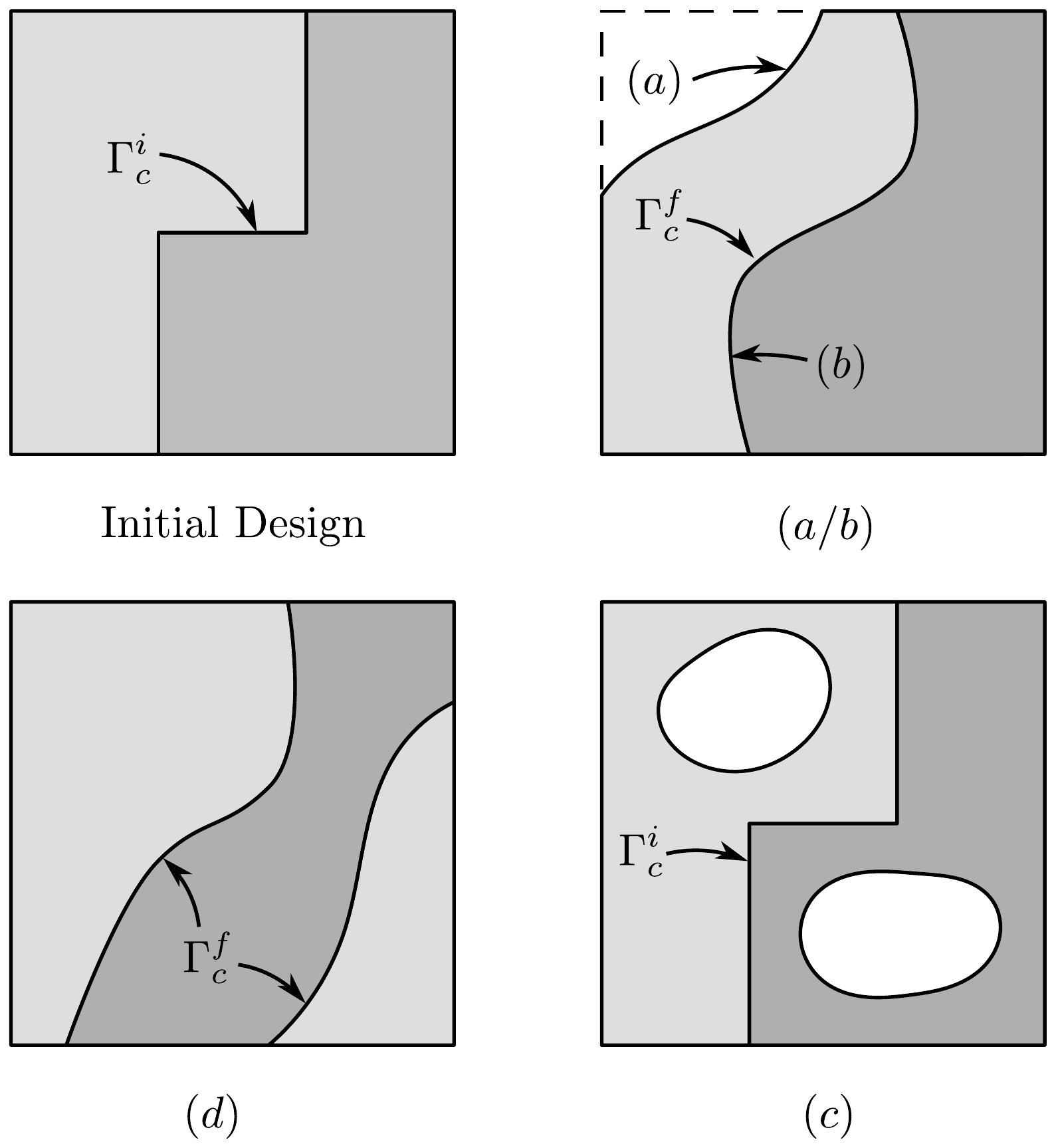}
\caption{Classifications of geometry control in contact optimization problems; $\Gamma_c^i$ and $\Gamma_c^f$ represent the initial and final contact interface geometry, respectively.}
\label{fig:topeq}
\end{center}
\end{figure}

Design optimization methods for contact related problems can be categorized by the type of geometry changes allowed during optimization. Figure \ref{fig:topeq} illustrates an initial design configuration comprised of two components in contact at the material interface, $\Gamma_c^i$, and four distinct options for geometry control. The shape of the external geometry excluding the interface can be altered (Fig.~\ref{fig:topeq}a), the shape of the interface geometry can be optimized (Fig.~\ref{fig:topeq}b), the design topology excluding the interface geometry can be varied (Fig.~\ref{fig:topeq}c), and the topology of the design including the interface geometry can be manipulated (Fig.~\ref{fig:topeq}d). While the proposed framework is flexible to use any of these geometry control options, this work studies option (d) for two-phase problems, and a combination of options (b) and (c) for three-phase problems.

In finite strain contact mechanics, the kinematic and constitutive nonlinearities which describe equilibrium are typically smooth and differentiable, whereas the interface conditions for contact and separation introduce a sharp discontinuity. Contact forces only act to prevent the interpenetration of bodies but vanish if the bodies separate. In addition to this sharp discontinuity, contact forces depend on surface orientation. For problems exhibiting large relative motion between components, coincident surface location is deformation dependent, complicating the solution of the physical response and evaluation of design sensitivities. These complexities pose interesting challenges for shape and topology optimization.

The optimization of contact related problems has received much interest within the engineering community. A wealth of literature exists for unilateral contact optimization problems; for review of advances prior to the turn of the century, the reader is referred to \cite{HKP:99}. In more recent studies, shape optimization excluding contact surface geometry (Fig.~\ref{fig:topeq}a) has been achieved using adaptive mesh refinement techniques for small strain \cite{HTK:01} and large strain problems \cite{HLDS:00}. While conformal meshing retains a sharp definition of the material interface, geometry updates afforded by mesh refinement are computationally expensive and are known to cause sensitivity inconsistencies due to changing discretization \cite{SMR:00}.

Departing from conformal mesh optimization methods, density methods, such as the Solid Isotropic Material with Penalization (SIMP) method, have become a popular alternative approach. Originally developed by \cite{Bendsoe:89} and \cite{RZB:92} for structural topology optimization, the SIMP method describes the geometry as a material distribution within the design domain. A fictitious porous material with density, $\rho$, is introduced to allow a continuous transition between two or more material phases. For more information on density methods and recent developments, the reader is referred to \cite{BS:03}, \cite{SM:13}, and \cite{DG:14}. The continuous density distribution within the design domain effectively smears the interface geometry, complicating the evaluation of design dependent surface loads. A common approach to circumvent this issue is to convert surface loads into volumetric body forces; see for example \cite{KC:01a}, \cite{SC:07}, and \cite{Yoon:10}. However, this method does not explicitly define the interface geometry and is ill-suited for modeling contact behavior. A second approach introduced by \cite{HO:00} is to apply surface loads by estimating the boundary via iso-volumetric density curves. This technique introduces approximation errors in interface position and orientation, degrading the reliability of the physical response prediction. A third approach, which has proven successful in the topology optimization of contact problems, is to provide an interface conforming mesh and optimize the surrounding material distribution. Analogous to Figure \ref{fig:topeq}d, topology changes have been afforded through density methods in small strain \cite{ARS:12,S:09} and large strain \cite{LLK:15}, excluding the contact surface from geometry control. This, however, severely restricts the optimal design solution space, as the functionality is often strongly correlated to the interface geometry.

Level set methods (LSM) provide a promising alternative approach to density methods as they retain a clear definition of the interface geometry. The interface is defined explicitly as the iso-contour of the Level Set Function (LSF) $\phi$ at a particular value, commonly $\phi=0$. For a review of recent developments of LSMs, the reader is referred to \cite{DML+:13}. The interface geometry is represented in the discretized mechanical model either via a body fitted mesh, an Ersatz material approach, or immersed boundary techniques.  The optimization of unilateral contact surface geometries (similar to Fig.~\ref{fig:topeq}b) have been achieved with LSM for small strain theory problems; see for example \cite{M:09}. Topology optimization including the material interface geometry has been achieved in a few small strain theory studies, namely for frictionless two-phase problems \cite{LM:15} and cohesive interface phenomena of multi-phase problems \cite{LLK:16}.  These two studies analyzed two dimensional problems and are comparable to option (d) and a combination of options (c) and (d) from Figure \ref{fig:topeq}, respectively.

Previous studies of optimization in which the contact interface is altered rely on either small strain kinematics or unilateral contact to reduce the complexity of sliding contact behavior. In this paper we expand the methods presented in \cite{LM:15} to the shape and topology optimization of bilateral contact problems with finite strain kinematics and large sliding contact. This marks a significant extension to the limits of accurate physical response prediction, which in turn grants access to a much broader scope of engineering problems. The proposed method allows altering both shape and topology of the contact interface. While geometric features can merge, no phase can be nucleated within volumes occupied by another phase. We use an explicit level set method to describe the interface geometry between two distinct material phases. In contrast to advancing the LSF by the Hamilton-Jacobi equation, explicit LSMs treat the parameters of the discretized LSF as explicit functions of the optimization variables \cite{DML+:13}.  This allows solving the resulting optimization problem by standard nonlinear programming algorithms.

We adopt an immersed boundary technique, specifically the eXtended Finite Element Method (XFEM), for predicting the mechanical response. The reader is referred to \cite{FB:06} and \cite{K:15} for an introduction and general overview of the XFEM. Modeling contact problems with the XFEM has shown great promise considering both friction and sliding contact. Assuming infinitesimal strains, sliding contact behavior has been modeled with the XFEM using penalty methods \cite{LB:08,MWL:12}, Lagrange multiplier methods \cite{GMM:07,AHL:12,LB:10}, and mortar methods \cite{GTT+:10}. Additionally, the XFEM has been leveraged to analyze problems in which relative sliding is significant. Large sliding bilateral contact behavior was considered using an augmented Lagrange method and surface-to-surface (STS) integration in small strain \cite{SGM+:13} and hybrid elements in large strain theory \cite{NGM+:09}. In large strain theory, penalty methods have proven successful for unilateral contact problems \cite{BP:12} and bilateral contact problems with node-to-surface (NTS) integration \cite{TM:11}. In this study we adopt a large strain theory stabilized Lagrange multiplier method similar to the approach of \cite{NGM+:09}. However, instead of using hybrid elements, contact equilibrium is enforced weakly through STS integration at the immersed boundary.

The remainder of this paper is organized as follows: in Section \ref{sec:opt}, we outline the formulation of the optimization problems considered in this study. In Section \ref{sec:fw}, we discuss the geometry model to describe the phase boundaries. In Section \ref{sec:phys}, the mechanical model of the contact problem is described. The XFEM formulation is summarized in Section \ref{sec:xfem}. Numerical implementation considerations are discussed in Section \ref{sec:num_impl}. In Section \ref{sec:ex}, we study the main characteristics of the proposed LSM-XFEM method with numerical examples. Insight gained from the numerical studies and areas for future research are summarized in Section \ref{sec:con}.

%===============================================================================

\section{Optimization Problem}\label{sec:opt}
\begin{figure}[t]
\begin{center}
\includegraphics[width=.8\linewidth]{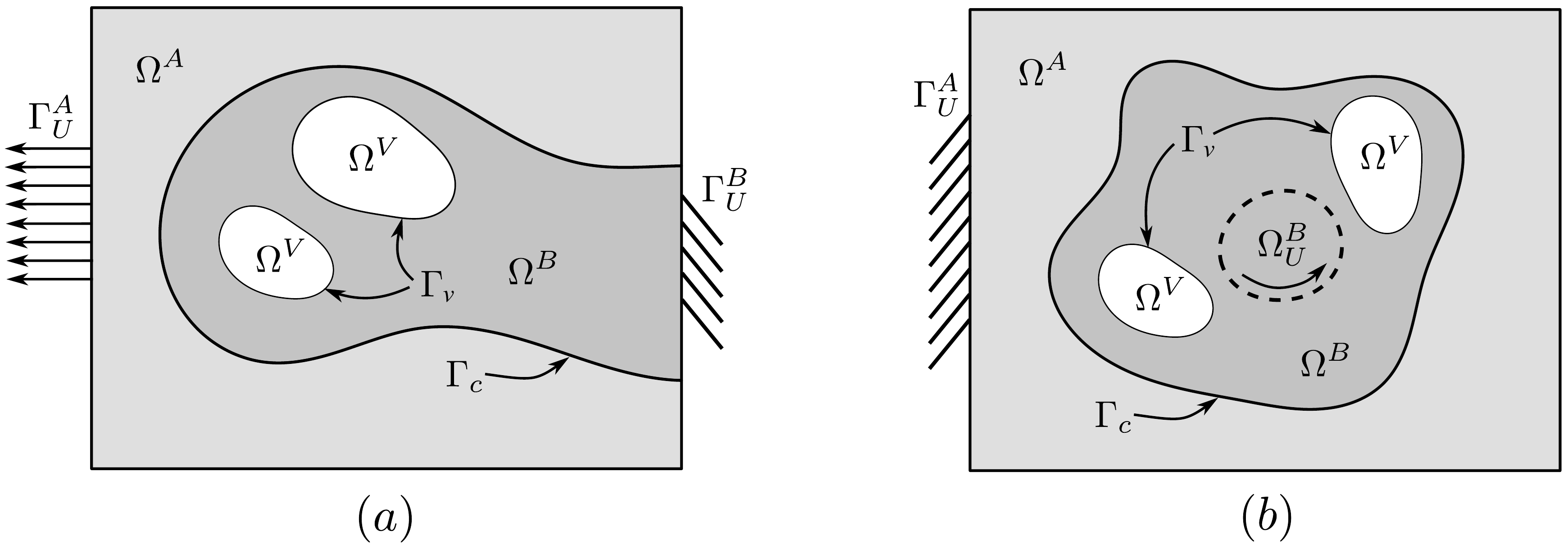}
\caption{Representative configurations of optimization problems pertinent to this study.}
\label{fig:opt_problem}
\end{center}
\end{figure}

\vspace{-2pt}
In this study we consider the interactions between two solid phases, A and B, with sliding, separable contact at the phase boundaries. For select optimization problems, a void phase, V, is introduced within solid phase B. The optimization problems presented in this paper can be illustrated by the representative configurations provided in Figure \ref{fig:opt_problem}. The design domain $\Omega_D$ is composed by three non-overlapping subdomains, $\Omega^{A}$, $\Omega^{B}$, and $\Omega^V$ such that $\Omega_D = \Omega^{A} \cup \Omega^{B} \cup \Omega^V$. The contact interface $\Gamma_C$ resides between the two solid phases such that $\Gamma_C= \Omega^{A} \cap \Omega^{B}$. The boundary between phase B and the void phase is denoted by $\Gamma_v= \Omega^{B} \cap \Omega^{V}$. To reduce interface complexities, such as triple junctions, the void subdomain, $\Omega^V$, resides within $\Omega^B$ such that $\Omega^V \cap \Omega^A = 0$. The approach for restricting the void phase to reside within phase B is discussed in Section \ref{sec:fw}.

While the proposed optimization method is applicable to a broad range of problems, we focus in this paper on two representative types of problems depicted in Figure \ref{fig:opt_problem} and studied in Section \ref{sec:ex}. The examples presented in Sections \ref{sec:anchor} and \ref{sec:snapfit} are analogous to Figure \ref{fig:opt_problem}(a), wherein the displacements in phase B are prescribed along the boundary $\Gamma^B_U$ and displacement controlled loading is applied at the boundary $\Gamma^A_U$. We seek to minimize an objective function related to the reaction load at $\Gamma_{U}^B$. The examples presented in Section \ref{sec:torqueprob} are analogous to Figure \ref{fig:opt_problem}(b), wherein the displacements in phase A are prescribed along the boundary $\Gamma^A_U$ and displacement controlled loading is applied within a subset of the domain occupied by phase B, $\Omega_U^B$. For problems considered here, the objective is to minimize some function related to the reaction load at boundary $\Gamma^A_U$.

The design problems of interest are defined by the following nonlinear program:
\begin{equation}
\label{eq:opt}
\begin{aligned}
                   & \min_{\mathbf{s}} q \left( \mathbf{s} \right), \\
    \text{s.t.} \  &  \frac{ V^B(\mathbf{s})}{V^B(\mathbf{s}) +  \  V^A(\mathbf{s})}  - c_v \leq  0 \\
                   & \mathbf{s} \in \mathbf{S} = \left\{ \mathbb{R}^{N_s} | s_{min} \leq s_i \leq s_{max}, \ i=1.... N_s  \right\} \ ,
\end{aligned}
\end{equation}
where $q$ denotes the scalar objective, $\mathbf{s}$ is the vector of optimization variables, and the number of optimization variables is $N_s$; the lower and upper bounds on the optimization variables are denoted by $s_{min}$ and $s_{max}$, respectively. For the scope of optimization problems studied in this paper, the objective function, $q$, is defined as:
\begin{equation}
q \left( \mathbf{s} \right) = c_u \frac{z \left(\mathbf{s},\hat{\mathbf{u}}(\mathbf{s})\right)}{z_0}+ c_p  \ \frac{P(\mathbf{s})}{P_0}
\end{equation}
where $z$ denotes the contribution of the mechanical response to the objective, $c_u$ is the associated weighting factor. To discourage the emergence of small geometric features, we introduce a perimeter penalty term, $P$, into the formulation of the objective function. The perimeter measure is the interface area of $\Gamma_c$ and $\Gamma_v$, and is computed as follows:
 \begin{equation}
 P = \int_{\Gamma_c \cup \Gamma_v}\mathrm{d}\Gamma.
 \end{equation}
The mechanical response contribution and the perimeter measure penalty are normalized by the initial measures, $z_0$ and $P_0$, respectively. While a perimeter penalty does not explicitly control the local shape and the feature size, it has been reported effective in regularizing structural optimization problems \cite{DML+:13}. In addition, for specific problems we constrain the ratio of volumes occupied by either solid, $V^A$ and $V^B$, to exclude trivial solutions. To provide control over the weighting of both the perimeter penalty and the volume inequality constraint, $c_p$ is the weight of the perimeter penalty, and $c_v$ controls the desired volume ratio between the two solids. While the proposed optimization framework allows considering other objectives and constraints, such as strain energy, displacement and stress measures, we found that the formulations of the optimization problem used here are well suited to illustrate the influence of the interface condition on the optimized design.

The dependency of the objective function and constraints on the optimization variables, $\mathbf{s}$, are defined by the framework described in Section \ref{sec:fw}. Note that the objective also depends on the structural response: $z(\mathbf{s}, \hat{\mathbf{u}})$, where $\hat{\mathbf{u}}$ denotes the vector of discretized state variables that are considered dependent variables of $\mathbf{s}$, i.e.~$\hat{\mathbf{u}}(\mathbf{s})$. The discretized state equations are described in Section \ref{sec:phys}. The optimization problem is solved by a nonlinear programming (NLP) method, and the design sensitivities are calculated by the adjoint method.

%===============================================================================

\section{Geometry Model}\label{sec:fw}
\vspace{-2pt}
\subsection{Two-Phase Problems}
The material layout of a two-phase problem is described by a LSF, $\phi(\mathbf{s},\mathbf{X})$, as follows:
\begin{equation}
\label{eq:lsm}
\begin{aligned}
    \phi(\mathbf{s},\mathbf{X})  & < 0,  && \forall \ \mathbf{X} \in \Omega^A \ , \\
    \phi(\mathbf{s},\mathbf{X})  & > 0,  && \forall \ \mathbf{X} \in \Omega^B \ , \\
    \phi(\mathbf{s},\mathbf{X})  & = 0,  && \forall \ \mathbf{X} \in \Gamma_c \ ,
\end{aligned}
\end{equation}
where $\mathbf{X}$ is the vector of spatial coordinates. Instead of updating the LSF by the solution of the Hamilton-Jacobi equation, as proposed by \cite{WWG:03} and \cite{AJT:04}, in this work the parameters of the discretized LSF are defined as explicit functions of the optimization variables.

For two-phase problems we follow the approach of \cite{KM:12} and discretize the design domain by finite elements and associate one optimization variable with each node, i.e.~$N_s = N_n$, where $N_n$ is the number of nodes. The level set value at the $i^{th}$ node is defined by the following linear filter:
\begin{equation}
\label{eq:levelsetFilter}
\phi_i = \left({\sum\limits^{N_n}_{j=1} w_{ij}}\right)^{-1} \ \sum\limits^{N_n}_{j=1} w_{ij}  s_j  \ ,
\end{equation}
with
\begin{equation}\label{eq:smooth_rad}
w_{ij}={max}\left(0, (r_f-|\mathbf{X}_i - \mathbf{X}_j|)\right) \ ,
\end{equation}
where $r_f$ is the filter radius, and $\mathbf{X}_j$ the position of the $j^{th}$ node. The level set filter \eqref{eq:levelsetFilter} widens the zone of influence of the optimization variables on the LSF and thus enhances the convergence of the optimization process \cite{KM:12}.

\subsection{Three-Phase Problems}
\label{sec:3phasegeom}

A popular approach to defining the spatial distribution of multiple materials with the LSM is through the superposition of multiple LSFs. Originally developed for digital image processing \cite{VC:02}, this method describes the layout of $2^m$ materials with $m$ LSFs. The individual phases are defined by the set of signs of the LSFs; the interfaces are described by one of the LSFs being zero. Also known as the `color' level sets method, this method has been reported useful in several multi-phase optimization studies \cite{DAP+:12,WW:05,YX:04a}.

In this study we take a similar approach by using two LSFs to distinguish three material phases; however, we limit the spatial arrangement of these phases as follows:
\begin{equation}\label{eq:geo_3phase}
\begin{aligned}
   \phi^1(\mathbf{s},\mathbf{X})& < 0,   &&                                      &&& \forall \ \mathbf{X} \in \Omega^A,   \\
   \phi^1(\mathbf{s},\mathbf{X})& > 0, \ && \phi^2(\mathbf{s},\mathbf{X}) > 0,   &&& \forall \ \mathbf{X} \in \Omega^B,   \\
   \phi^1(\mathbf{s},\mathbf{X})& > 0, \ && \phi^2(\mathbf{s},\mathbf{X}) < 0,   &&& \forall \ \mathbf{X} \in \Omega^{V}, \\
   \phi^1(\mathbf{s},\mathbf{X})& = 0,   &&                                      &&& \forall \ \mathbf{X} \in \Gamma_c,   \\
   \phi^2(\mathbf{s},\mathbf{X})& = 0,   &&                                      &&& \forall \ \mathbf{X} \in \Gamma_v.
\end{aligned}
\end{equation}
This conditional treatment of the LSFs admits the definition of a third phase; however, $\Omega^{V}$ is restricted to reside within phase B through \eqref{eq:geo_3phase}. For this work, the LSFs $\phi^1$ and $\phi^2$ are parameterized to describe a set of geometric primitives, such as circles or rectangles. The optimization variables define the location and the dimensions of the primitives. To avoid the emergence of triple junctions where $\phi^1=\phi^2 = 0$, the limits $s_{min}$ and $s_{max}$ are carefully selected to ensure the geometric primitives in $\phi^1$ do not intersect the geometric primitives in $\phi^2$. This approach is used in two examples presented in Section \ref{sec:ex}.

%===============================================================================

\section{Physics Model}\label{sec:phys}
\vspace{-2pt}
Static equilibrium of phases $\Omega^A$ and $\Omega^B$ within the design domain is satisfied by the balance of linear momentum referred to the reference configuration $\Omega^p_0$ for $p=A,B$:
\begin{equation}
\label{eq:bom}
 \nabla \cdot \left(\mathbf{F}^p \ \mathbf{S}^p \right) +  \mathbf{b}^p_0 = \mathbf{0} \quad \textrm{in} \  \Omega^{p}_0 \ ,
\end{equation}
subject to the Dirichlet boundary conditions:
\begin{equation}
  \mathbf{u}^p = \mathbf{U}^p \quad \textrm{on} \ \Gamma^p_U \ ,
\end{equation}
where $\mathbf{u}^p$ is the displacement vector, $\mathbf{F}^p$ is the deformation gradient tensor, $\mathbf{S}^p$ is the second Piola-Kirchhoff stress tensor, $\mathbf{b}^p_0$ is the reference configuration body force vector, and $\mathbf{U}^p$ is the vector of prescribed displacements at the boundary $\Gamma^p_U$. We assume a hyper-elastic neo-Hookean material behavior and a nonlinear kinematic relationship:
\begin{equation}
\mathbf{S}^p = \lambda^p \ \ln \left({ \det \mathbf{F}^p }\right) \ {\mathbf{C}^p}^{-1} + \mu^p \left({ \mathbf{I}-{\mathbf{C}^p}^{-1} }\right),
\end{equation}
with
\begin{equation}
       \mathbf{C}^p = {\mathbf{F}^p}^T \ \mathbf{F}^p,
\qquad \mathbf{F}^p = \frac{\partial \mathbf{x}^p}{\partial \mathbf{X}^p},
\qquad \mathbf{x}^p = \mathbf{u}^p + \mathbf{X}^p,
\end{equation}
where $\lambda^p$ and $\mu^p$ represent the material Lam\'e parameters, $\mathbf{C}^p$ is the right Cauchy-Green tensor, $\mathbf{I}$ is the identity matrix, $\mathbf{x}^p$ is the current position, and $\mathbf{X}^p$ is the reference position of phase $p=A,B$.

\begin{figure}[t]
\begin{center}
\includegraphics[width=.6\linewidth]{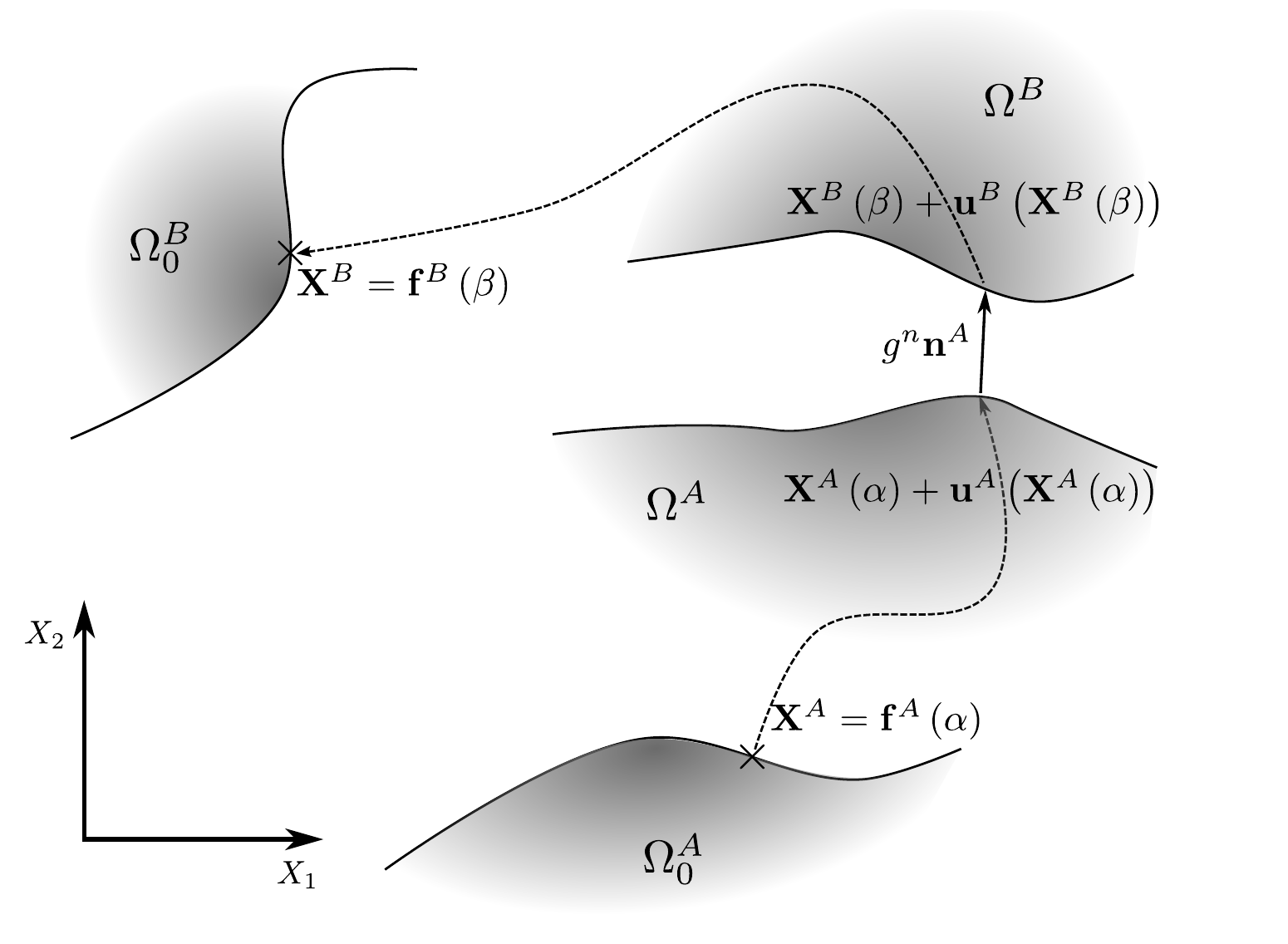}
\caption{Parametric representation of surfaces belonging to materials A and B.}
\label{fig:pmapping}
\end{center}
\end{figure}

In the presence of large relative motion between surfaces, the dependence of coincident location along the interface on the displacements of either body needs to be accounted for. To this end, the surfaces of both structural phases are mapped to parametric space. This parametrization simplifies the definition of coincident surface location by describing surface position $\mathbf{X}^p$ and subsequently the displacements $\mathbf{u}^p$ in a reduced dimensional space. The surfaces of phase $A$ and $B$ are parameterized by some functions $\mathbf{f}^A$ and $\mathbf{f}^B$ of the surface parameters $\alpha$ and $\beta$ respectively, as illustrated in Figure \ref{fig:pmapping}.
The set of parametric functions, $\mathbf{f}^A$ and $\mathbf{f}^B$, used in this paper are directly related to the method of discretization, and are detailed in Section \ref{sec:xfem}.

To provide a continuous representation of coincident surface position, both surface parametrization schemes are coupled through the following relationship:
\begin{equation}\label{eq:par_map}
\mathbf{X}^A\left(\alpha\right) + \mathbf{u}^A\left(\mathbf{X}^A\left(\alpha\right)\right) + g_n \mathbf{n}^A\left(\mathbf{X}^A,\mathbf{u}^A\right) - \mathbf{X}^B\left(\beta\right) - \mathbf{u}^B\left(\mathbf{X}^B\left(\beta\right)\right) = \mathbf{0} \ ,
\end{equation}
where $g_n$ is the magnitude of the gap between both surfaces in the direction of the deformed configuration surface normal $\mathbf{n}^A$; see Figure \ref{fig:pmapping}. Utilizing a master-slave approach, surface position $\beta$ and the scalar normal gap $g_n$ are defined through \eqref{eq:par_map} for any given surface position $\alpha$. Qualitatively, this expression states that for any given position along surface $\mathbf{x}^A$ in the current configuration, the coincident position along surface $\mathbf{x}^B$ can be found by a projection in the direction of the deformed surface normal, $\mathbf{n}^A$.

Along either surface in the reference configuration, the following non-penetration conditions apply:
\begin{equation}
\label{eq:non-pen}
g^p_0 \ \lambda^p_0 = 0, \ g^p_0 \geq 0, \ \lambda^p_0 \leq 0,
\end{equation}
with
\begin{equation}
\label{eq:contactpres}
\lambda^p_0 =  \left(\mathbf{n}^{p}_0\right)^T \mathbf{S}^p \ \mathbf{n}^{p}_0 ,
\end{equation}
\begin{equation}
\label{eq:mapgap}
g^p_0 = g^n j^p ,
\end{equation}
\begin{equation}
j^p = \det\left(\mathbf{F}^p\right) \| {\mathbf{F}^p}^{-T} \mathbf{n}^{p}_0 \|,
\end{equation}
where the $g^p_0$ is the normal gap between the bodies pulled back to the reference configuration of material $p$, $g^n$ is the normal gap between bodies in the deformed configuration, and $\lambda^p_0$ is the surface traction in normal direction in the reference configuration. The Jacobian of the surface area, $j^p$, is derived from Nanson's formula, as outlined in \cite{Wriggers:02}. As the bodies cannot interpenetrate, the gap must be greater than or equal to zero. The surface traction is negative when bodies are in contact, but vanishes as they separate. Thus, $\lambda^p_0$ serves as the Lagrange multiplier of the non-penetration condition. Considering that in the deformed configuration the surface pressures are identical,
\begin{equation}
\lambda_B = \lambda_A \quad \textrm{and thus} \quad \lambda^B_0{j^B}^{-1} = \lambda^A_0 {j^A}^{-1},
\end{equation}
we express the surface pressure with just $\lambda^A_0$, residing within the master reference configuration, $\Gamma^A_{c,0}$. In particular, $\Gamma^A_{c,0}$ is the undeformed contact surface of phase A. To simplify notation, we drop the superscript and define the Lagrange multiplier as $\lambda_0 \equiv \lambda^A_0$. Provided that in three-phase problems we do not allow the void material interface to directly connect to the contact interface $\Gamma_c$ , i.e.~$\phi^1 = \phi^2 = 0$, both material surfaces are coincident in the reference configuration and the initial gap, $g_n$, is zero.

The XFEM discretization of the contact problem is based upon the following stabilized weak form of the governing equations:
\begin{multline}
\label{eq:weakform}
  \sum_{p=A,B}{\int_{\Omega^{p}_0}{ \mathbf{F}(\bm{\nu}^p) :  \left( {\mathbf{F}^p \mathbf{S}^p } \right) \ d \Omega}}
- \sum_{p=A,B}{\int_{\Omega^{p}_0}{ \bm{\nu}^p \cdot \mathbf{b}_0^p \ d \Omega}}  \\
- \sum_{p=A,B}{\int_{\Gamma^{p}_{N,0}}{ \bm{\nu}^p \cdot \mathbf{T}_0^p \ d \Gamma}}
- \int_{\Gamma^{A}_{c,0}}{ \delta g^A_0 \ \lambda_0 \ d \Gamma} + r^G= 0 \ ,
\end{multline}
where $\bm{\nu}^p$ is an admissible test function, $\mathbf{T}_0^p$ is a prescribed traction at the external boundary $\Gamma^{p}_{N,0}$, $\delta g^A_0$ is the variation of the normal gap pulled back to the undeformed surface of phase A, and $r^G$ is a stabilization term discussed in Section \ref{sec:ghost_dyn}. Similar to the augmented Lagrange formulation presented by \cite{Wriggers:02}, the Lagrange multiplier is governed by the following constraint equation:
\begin{equation}
\label{eq:lagmultgov}
\int_{\Gamma^{A}_{c,0}} { {\mu} \ \left( \lambda_0 - \tilde{\lambda_0} - \gamma \ g^A_0 \right) \ d \Gamma} = 0,
\end{equation}
with
\begin{equation}
\label{eq:lambdaprx}
\tilde{\lambda_0} = \kappa^A {\mathbf{n}^A_0}^T \mathbf{S}^A \mathbf{n}^A_0 + \kappa^B {\mathbf{n}^B_0}^T \mathbf{S}^B \mathbf{n}^B_0 {j^B}^{-1} j^A
\end{equation}
where ${\mu}$ is a test function for the non-penetration condition, $\tilde{\lambda_0}$ is a weighted average of the surface traction in the normal direction and $\kappa^p$ are weighting factors such that $\kappa^A+\kappa^B=1$. In our experience, the penalty factor $\gamma$ discourages penetration during the early stages of convergence, but becomes insignificant as equilibrium is achieved. The formulations for $\kappa^p$ and $\gamma$ are related to discretization, and are provided in Section \ref{sec:xfem}. The constraint equation \eqref{eq:lagmultgov} and contact contributions to the weak form of the equilibrium equations \eqref{eq:weakform} are integrated over $\Gamma^A_{c,0}$, and an active set strategy is used to handle the inequality constraint regarding surface separation.

%===============================================================================

\section{XFEM Discretization}\label{sec:xfem}
\vspace{-2pt}
The XFEM provides an elegant approach to discretize the weak form of partial differential equations where the geometry is described by a LSF. The space of the test and trial functions are augmented by enrichment functions to capture weak and strong discontinuities within elements intersected by the zero level set iso-contour. For the problems considered in this paper, the displacement field is discontinuous across the contact interface $\Gamma_c$. Therefore, a Heaviside enrichment is exclusively used in this work. We approximate the displacements $u_i(\mathbf{X})$  for two phase problems as follows:
\begin{equation}
u_i(\mathbf{X}) = \sum_{m=1}^M \left( {   H(-\phi(\mathbf{X})) \ \sum^{N_e}_{k} \ N_k(\mathbf{X}) \ \delta^{A,k}_{mq} \ u^{A,m}_{i,k} } { \ +\  H( \phi(\mathbf{X})) \ \sum^{N_e}_{k} \ N_k(\mathbf{X}) \ \delta^{B,k}_{mp} \ u^{B,m}_{i,k} } \right) \ ,
\end{equation}
with $H$ being the Heaviside function:
\begin{equation}\label{eq:heaviside}
    H(\phi) =
    \begin{cases}
        1 & \text{if } \phi > 0, \\
        0 & \text{if } \phi \leq 0.
    \end{cases}
\end{equation}
The number of enrichment levels is denoted by $M$, $N_e$ is the number of nodes in the element, $N_i(\mathbf{X})$ are the shape functions, $u^{p,m}_{i,k}$ is the degree of freedom of enrichment level $m$ at node $k$ interpolating the displacement $u_i$ in phase $p$. The Heaviside function turns on/off two sets of shape functions associated with the phases A and B. For each phase, multiple enrichment levels, i.e.~sets of shape functions, might be necessary to interpolate the displacements in multiple, physically disconnected regions without introducing spurious coupling and load transfer between disconnected regions of the same phase. The Kronecker delta $\delta^{p,k}_{mq}$ selects the active enrichment level $q$ for node $k$ such that the displacements at a point $\mathbf{X}$ are interpolated by only one set of degrees of freedom defined at node $k$, satisfying the partition of unity principle. For further description see \cite{MM:13}, \cite{TAY:03}, and \cite{TYHT:11}. In elements not intersected by the zero level set iso-contour, the displacement field is approximated by a standard finite element interpolation. The enrichment level for these elements is chosen to maintain a continuous displacement field across element boundaries.

For three-phase problems, we approximate the displacements $u_i(\mathbf{X})$ in phases A and B as:
\begin{equation}
u_i(\mathbf{X}) = \sum_{m=1}^M \left( {   H(-\phi^1(\mathbf{X})) \ \sum^{N_e}_{k} \ N_k(\mathbf{X}) \ \delta^{A,k}_{mq} \ u^{A,m}_{i,k} } { \ +\  H(\phi^2(\mathbf{X}))H( \phi^1(\mathbf{X})) \ \sum^{N_e}_{k} \ N_k(\mathbf{X}) \ \delta^{B,k}_{mp} \ u^{B,m}_{i,k} } \right) \ .
\end{equation}
The Heaviside function applied to the LSF $\phi^2$ serves to turn off the approximation in the void phase. Aside from this deviation, the displacement field is approximated in the same manner as for two phase problems. For both two and three phase problems, the intersected elements in phases A and B are triangulated for integration purposes.

%===============================================================================

\section{Numerical Implementation}\label{sec:num_impl}
\vspace{-2pt}

%===============================================================================

\subsection{Contact Equilibrium Contributions}

\begin{figure}[t]
\begin{center}
\includegraphics[width=0.9\linewidth]{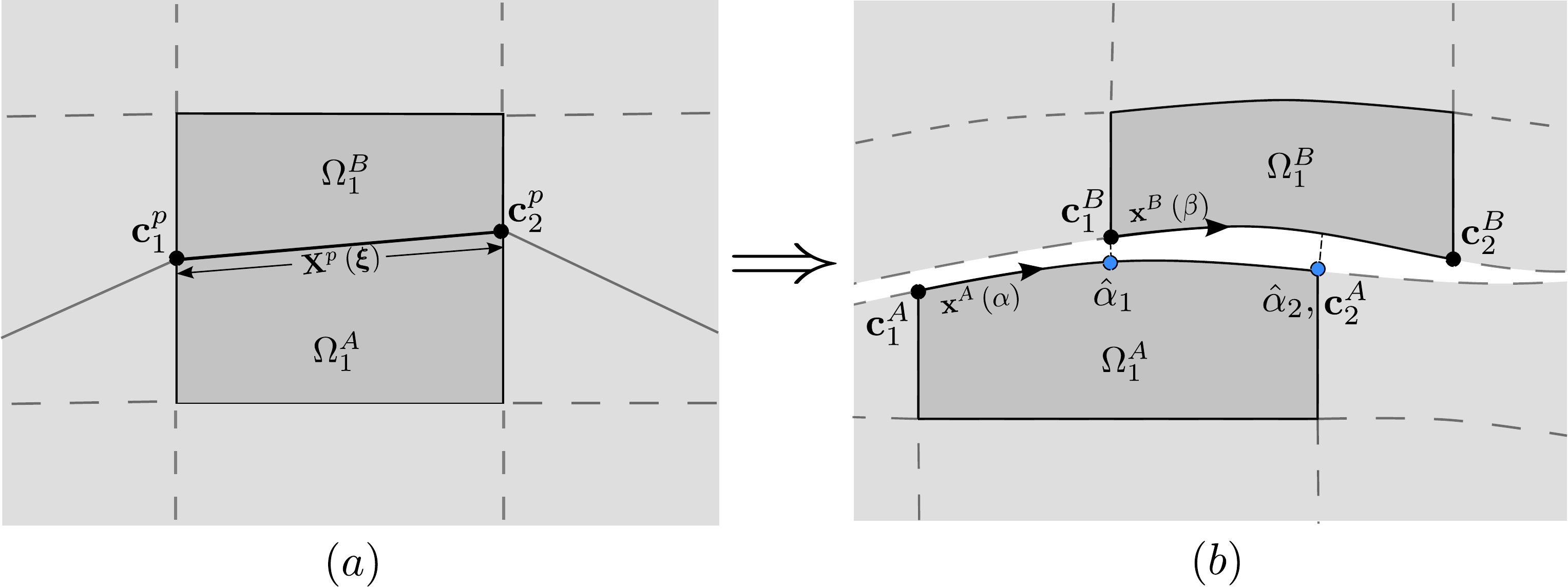}
\caption{Element immersed surface parametrization in the undeformed $(a)$} and current $(b)$ configuration.
\label{fig:cdofs}
\end{center}
\end{figure}

The XFEM retains a piece-wise continuous definition of the interface geometry subject to the chosen method of LSF interpolation. For STS integration of the contact contributions to the equilibrium equation \eqref{eq:weakform}, coincident locations at the contact interface must be identified as illustrated in Figure \ref{fig:cdofs}. Locations $\mathbf{c}_1^p$ and $\mathbf{c}_2^p$ correspond to the interface boundaries for a particular element of phase $p$, while $\hat{\alpha}_1$ and $\hat{\alpha}_2$ are the limits of integration for this particular element pair. Provided that the coincident surface location is governed by \eqref{eq:par_map}, element integration limits are deformation dependent. To recover a fully consistent tangent stiffness, which is essential to the accuracy of the adjoint sensitivity analysis, these integration limit dependencies on the solution must be accounted for.

If the integration limit coincides with the element boundary of phase A, i.e. $\mathbf{x}^A\left(\hat{\alpha}_i\right)=\mathbf{c}_i^A$, it is solution independent. However, if the integration limit $\hat{\alpha}_i$ does not coincide with the elemental boundary $\mathbf{c}_i^A$, as is the case for $\hat{\alpha}_1$ in Figure \ref{fig:cdofs}, its position depends on the projection of $\mathbf{c}_1^B$  onto the phase A surface in the deformed configuration. The integration limit $\hat{\alpha}_1$ and its dependencies on the displacement field are defined through \eqref{eq:par_map}.

In this paper, the zero level set iso-contour is interpolated linearly within an intersected element and the position of a point on the intersection is parameterized by:
\begin{equation}
\mathbf{X}^p\left(\xi\right) = \left(1-\xi\right)\mathbf{c}_1^p + \xi \mathbf{c}_2^p \ ,
\end{equation}
where the local coordinate $\xi$ corresponds to either $\alpha$ or $\beta$ for phase A or B respectively. The test and trial functions for the normal surface traction, ${\mu}$ and $\lambda_0$, are piecewise linear for each STS element pair, but not necessarily continuous across element boundaries. Thus, the associated degree of freedom can be computed for each STS element pair and condensed from the global system of equations.

Following the work of \cite{AHD:14}, the weighting factors $\kappa^p$ for computing the average normal traction in \eqref{eq:lambdaprx} depend on the elemental intersection configuration as follows:
\begin{equation}\label{eq:aweight}
\begin{aligned}
\kappa^A &= \frac{|\Omega|^A / E^A}{|\Omega|^A / E^A + |\Omega|^B / E^B} \ , \\
\kappa^B &= \frac{|\Omega|^B / E^B}{|\Omega|^A / E^A + |\Omega|^B / E^B} \ ,
\end{aligned}
\end{equation}
where $|\Omega|^p$ denotes the elemental volume occupied by phase $p=A,B$ and $E^p$ is the Young's modulus of phase $p=A,B$. The penalty factor in \eqref{eq:lagmultgov} depends on the element size $h$ and is set to:
\begin{equation}\label{eq:uniteq}
\gamma= \frac{E^A + E^B}{h} \ .
\end{equation}

%===============================================================================

\subsection{Stabilization}\label{sec:ghost_dyn}

During the design optimization process, the interface of the embedded geometry may produce intersection configurations where certain degrees of freedom interpolate to very small subdomains. This causes an ill-conditioning of the mechanical model system, which may impede the convergence of the nonlinear problem. In the context of contact problems, the vanishing zone of influence of a particular degree-of-freedom may also result in artificially high stress approximations in localized regions near the interface, leading to the erroneous evaluation of the contact pressure. To mitigate ill-conditioning of the system and poor structural response prediction at the interface, we apply a face-oriented ghost-penalty formulation. Similar to the stabilization method for diffusion problems presented by \cite{BH:12}, we penalize the jump in stress across element borders, and define the stabilization term $r^G$ in \eqref{eq:weakform} as:
\begin{equation}
r^G = \sum_{p=A,B}\int_{\Gamma_e^0} \gamma^{G}\left\llbracket \frac{\partial \bm{\nu}^p}{\partial \mathbf{X}} \right\rrbracket \mathbf{n}^0_e \left\llbracket \mathbf{S}^p \right\rrbracket \mathbf{n}^0_e\mathrm{d}\Gamma ,
\end{equation}
where $\Gamma_e^0$ is the reference configuration boundary of intersected elements, $\gamma^G$ is a penalty parameter, $\bm{\nu}$ is an admissible test function, $\mathbf{n}^0_e$ is the reference configuration surface normal of the element boundary, and the jump operator,
\begin{equation}
\left\llbracket \zeta \right\rrbracket = \zeta\arrowvert_{\Omega^1_e} - \zeta\arrowvert_{\Omega^2_e},
\end{equation}
computes the difference of a particular quantity across the facet between two adjacent elements, $\Omega_e^1$ and $\Omega_e^2$.  The penalty parameter, $\gamma^G$, is defined as:
\begin{equation}
\gamma^G = \epsilon \ h
\end{equation}
where $\epsilon$ is a problem-specific scaling factor and $h$ is the element side length. The jump in stress is penalized across the entire element border, irrespective of where the interface intersects it. Face-oriented ghost penalization has been reported as being beneficial to various fluid flow related problems, including fluid-solid interactions \cite{BFH:06}, high Reynolds number flows \cite{SW:14}, and incompressible flows \cite{SRG+:14}.

%===============================================================================

\subsection{Dynamic Relaxation}

Contact problems often experience moments of neutral equilibrium, and can exhibit snap-through behavior. In this work the discretized mechanical model is solved using a Newton-Raphson iterative procedure, which may suffer from convergence difficulties in such scenarios. To mitigate these issues, we use a Levenberg-Marquardt \cite{M:78} type method for dynamic relaxation. Initially developed to solve non-linear least square problems, this algorithm has also been reported useful in reducing analysis instabilities caused by element distortion in compliant mechanism optimization problems \cite{K:09}. We adopt a similar approach by modifying the tangent stiffness matrix:
\begin{equation}
\tilde{\mathbf{J}} = \mathbf{J} + \beta\ \mathrm{diag}\left(\mathbf{J}\right) \ ,
\end{equation}
where $\mathbf{J}$ is the original tangent stiffness matrix, $\beta$ is the damping parameter, and $\tilde{\mathbf{J}}$ is the modified tangent stiffness matrix. The damping parameter is given an initial value, and adaptively increased or decreased by a factor of 10 depending on whether the satisfaction of equilibrium improves or deteriorates throughout the iterative solution procedure. The initial value of the damping parameter is $\beta = 0.01$ for the problems presented in this paper.

%===============================================================================

\section{Numerical Examples}\label{sec:ex}
\vspace{-2pt}
To demonstrate the accuracy of the proposed framework, we first verify both the physical response prediction and adjoint method. Subsequently, a comparative study explores optimal design improvements for large strain theory over small strain theory. Finally, problems with objectives that characterize the structural response during a quasi-static loading process illustrate the main characteristics of the proposed framework.

For all examples, we assume plane strain conditions and a quasi-static response. Finite strain kinematics and large sliding contact are used unless specified otherwise. The mechanical model is discretized with bilinear Quad-4 elements using the framework described in Section \ref{sec:xfem}. The nonlinear contact problems are solved by Newton's method with dynamic relaxation, using an active set strategy for the contact conditions. A drop of the residual of $10^{-6}$ relative to the initial residual is required, unless stated otherwise. Loads are applied incrementally, and a direct solver is used for the linearized sub-problems.

The parameter optimization problems are solved by the Globally Convergent Method of Moving Asymptotes (GCMMA) of \cite{Svanberg:02}. The parameters for the initial, lower, and upper asymptote adaptation are set to $0.5$, $0.7$, and $1.2$, respectively. The relative step size, $\Delta s$, is given with each example. The design sensitivities are computed with the adjoint method. The reader is referred to \cite{AK:96}, \cite{KPC:01} and \cite{KYC:02} for an in depth discussion of sensitivity analysis for contact problems. In this work, the partial derivatives of the state equations and objective function with respect to the state variables are evaluated using analytically differentiated formulations. The partial derivatives of the objective, constraints, and element residuals with respect to the optimization variables are calculated by a finite difference scheme, which is computationally inexpensive as only intersected elements need to be considered.

%===============================================================================

\subsection{Mechanical Model Verification}
\begin{figure}[b]
\begin{minipage}[t]{0.48\textwidth}
\centering
\vspace{0pt}
\includegraphics[width=.7\linewidth]{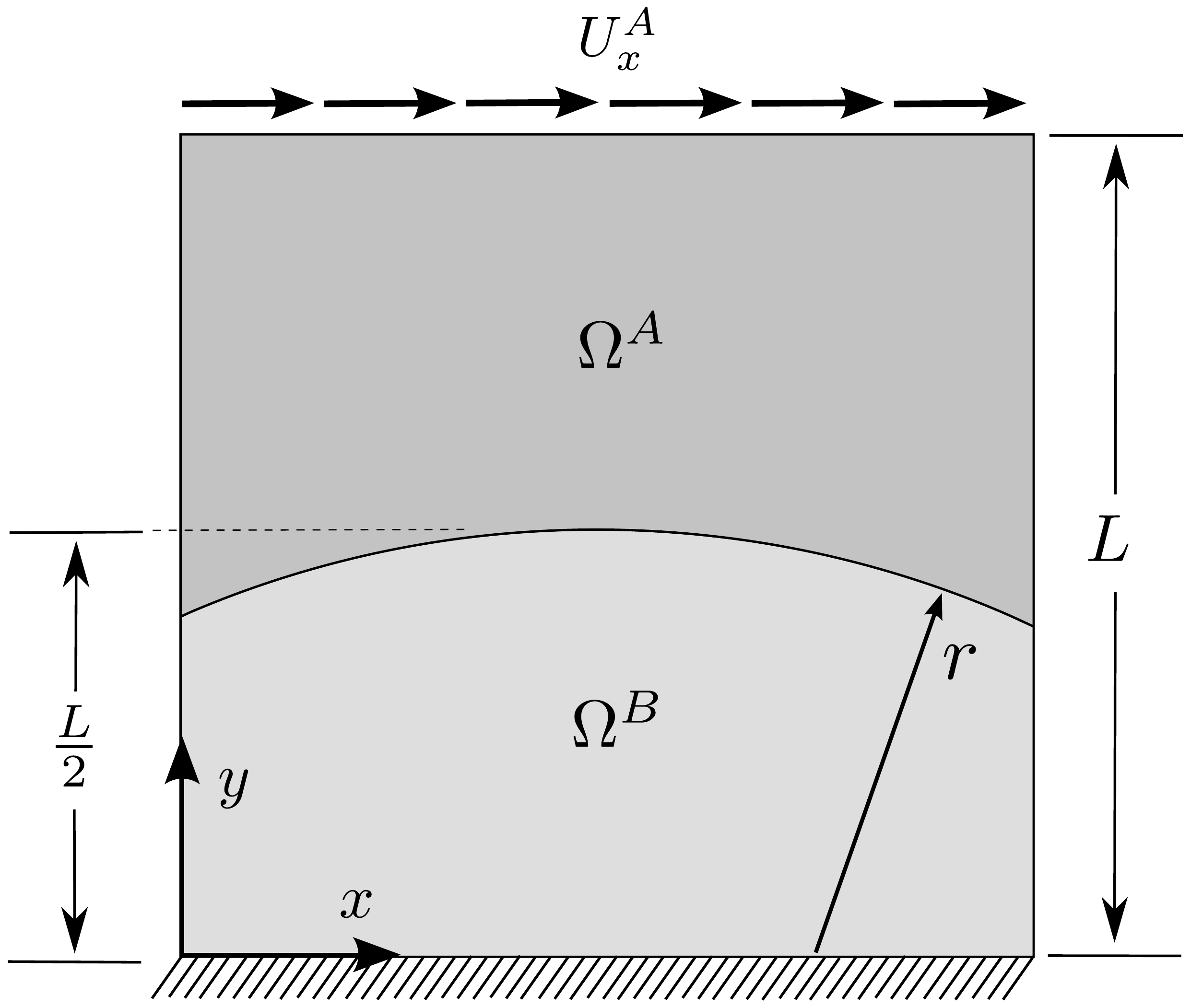}
\caption{Large strain, frictionless contact benchmark setup.}
\label{fig:bench-ls-fric-setup2}
\end{minipage}\hfill
\begin{minipage}[t]{0.48\textwidth}
\centering
\vspace{0pt}
\captionof{table}{Benchmark mechanical model parameters.}
\label{tab:bench-ls-fric-setup2}
\begin{tabular}{ll}
\hline
Description & Parameter \\ \hline
domain length & $L = 1.0$ m \\
interface radius & $r = 1.2$ m \\
Young's modulus & $E^A = 10$ MPa \\
Young's modulus & $E^B = 10$ MPa \\
Poisson's ratio & $\nu^A = 0.3$ \\
Poisson's ratio & $\nu^B = 0.3$ \\
applied displacement & $U_x^A = 0.5$ m \\
\hline
\end{tabular}
\end{minipage}
\end{figure}

To verify the accuracy of the XFEM mechanical model, a benchmark example is studied and compared to results produced by Abaqus\textsuperscript{\textregistered} using a conformal mesh. We consider a square domain which is composed of two non-overlapping subdomains $\Omega^A$ and $\Omega^B$; see Figure \ref{fig:bench-ls-fric-setup2}. The contact interface $\Gamma_c = \Omega^A \cap \Omega^B$ is defined by an arc of radius $r$. The volumes occupied by either phase, $\Omega^A$ and $\Omega^B$, are modeled by neo-Hookean materials of the same properties. Displacements at the top edge of $\Omega_A$ are prescribed and incrementally increased in $50$ load steps to a maximum value of $U_x^A = 0.5$ and $U_y^A = 0.0$. The bottom edge of $\Omega_B$ is fixed. Dimensions and material properties for the model are presented in Table \ref{tab:bench-ls-fric-setup2}.

\begin{figure}[t]
\begin{minipage}[b]{0.54\textwidth}
\centering
\vspace{0pt}
\includegraphics[width=1\linewidth]{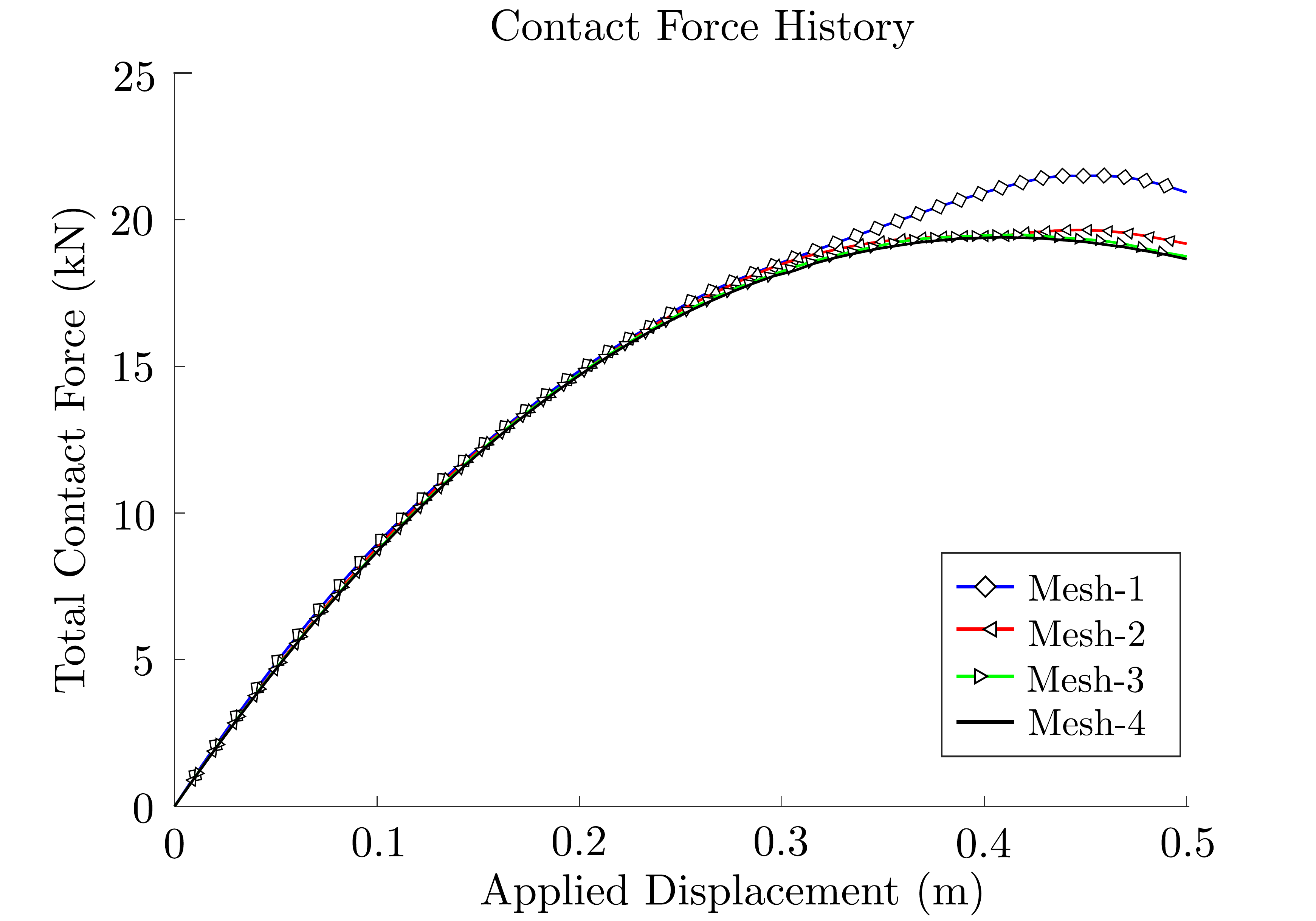}
\caption{Total surface force as a function of applied displacement.}
\label{fig:bench-ls-meshr_force2}
\end{minipage}\hfill
\begin{minipage}[b]{0.42\textwidth}
\vfill
\centering
\vspace{0pt}
\includegraphics[width=1\linewidth]{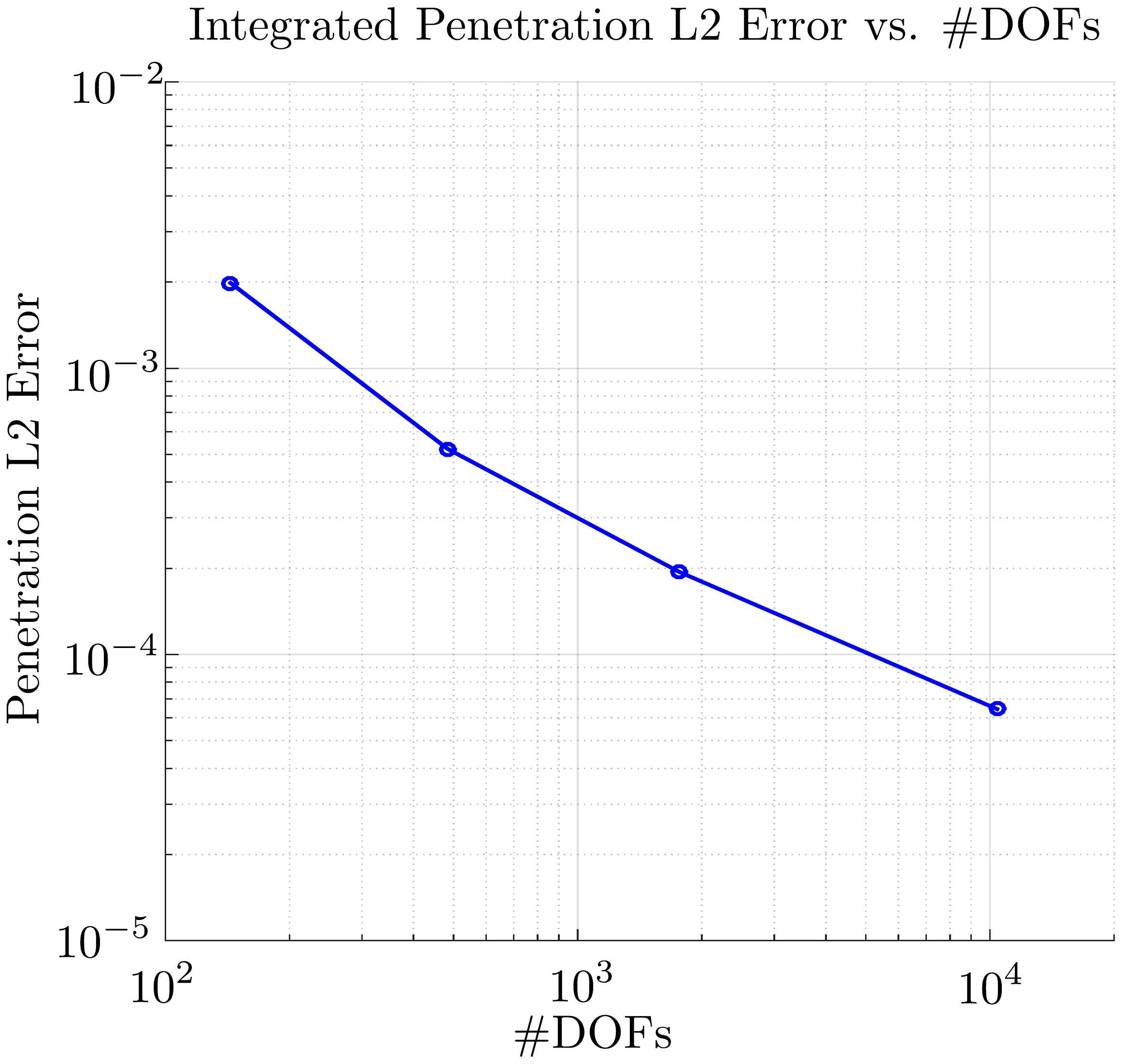}
\caption{Integrated $L_2$ error with mesh refinement.}
\label{fig:bench-ls-meshr_err2}
\end{minipage}
\end{figure}

To examine the convergence behavior of the mechanical model, the problem is analyzed for four different mesh sizes: Mesh 1 consists of $5\times5$ elements, Mesh 2 contains $11\times11$ elements, Mesh 3 has $21\times21$ elements, and Mesh 4 consists of $51\times51$ elements. The coupled parametric representation of coincident surface location facilitates a quadratic convergence, requiring on average 6 Newton iterations to converge to a tolerance criterion of $1\times10^{-9}$. Surface contact forces and normal penetration errors are extracted at the interface. Figure \ref{fig:bench-ls-meshr_force2} demonstrates the total contact force convergence with mesh refinement.

Additionally, Figure \ref{fig:bench-ls-meshr_err2} illustrates the integrated penetration $L_2$ error for each mesh used. The $L_2$ error is determined as follows:
 \begin{equation}
 \mathrm{L_2 \ error} = \sqrt{\frac{\int_{\lambda_0 < 0} \int_{\Gamma_c}\  g_n^2 \ \mathrm{d}\Gamma\mathrm{d}t}{\int_{\Gamma_c}\mathrm{d}\Gamma}} \ ,
 \end{equation}
where $t$ is a pseudo-time which describes the loading process, and $g_n$ is the normal gap between surfaces in contact. The normal gap, $g_n$, is only integrated across contact element pairs when they are in an active state of contact, i.e. $\lambda_0 < 0$.  The curved interface is described by a linearly interpolated LSF. This approach leads to a segmented interface that may yield poor response predictions at low levels of mesh discretization. As the mesh is refined, the force profile converges; see Figure \ref{fig:bench-ls-meshr_force2}. Considering the non-penetration condition is enforced weakly at the interface, surface penetration error diminishes with mesh refinement; see Figure \ref{fig:bench-ls-meshr_err2}.
\begin{figure}[t]
\begin{center}
\includegraphics[width=1\linewidth]{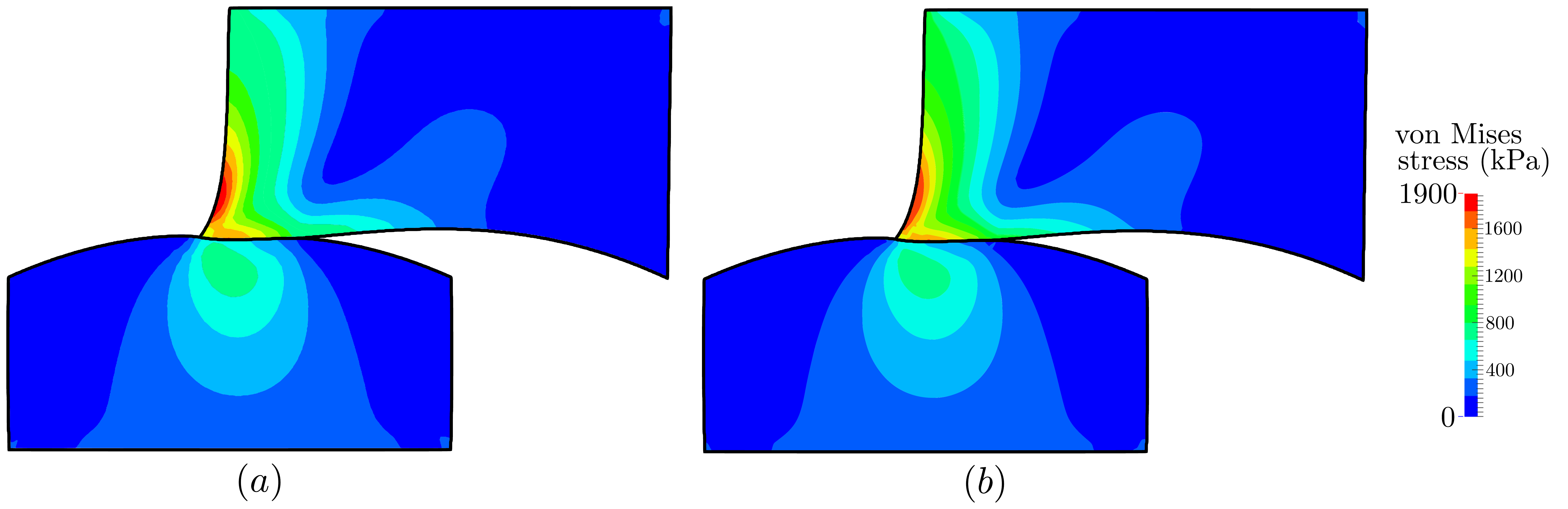}
\caption{Comparison of von Mises stress for $(a)$ Abaqus\textsuperscript{\textregistered} and $(b)$ current implementation.}
\label{fig:bench-ls-compare-vm2}
\end{center}
\end{figure}
\begin{figure}[t]
\begin{center}
\includegraphics[width=.55\linewidth]{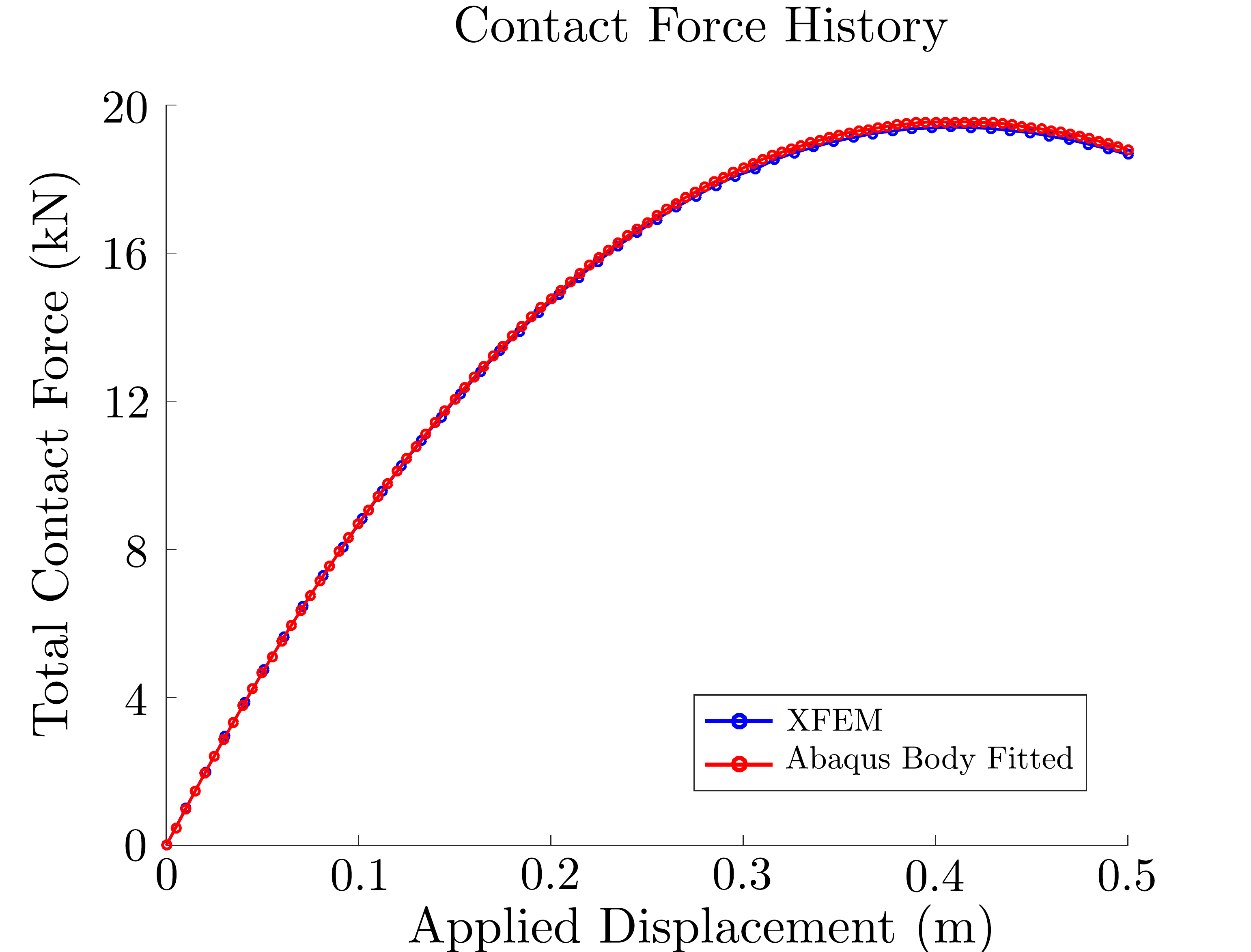}
\caption{Comparison of total contact force as a function of applied displacement.}
\label{fig:bench-ls-compare-cf2}
\end{center}
\end{figure}

For verification purposes, this benchmark problem is modeled in Abaqus\textsuperscript{\textregistered} using a body-fitted mesh with 50$\times$50 Quad-4 plane strain elements. Interface conditions in Abaqus\textsuperscript{\textregistered} are enforced using a surface-to-surface, frictionless augmented Lagrange formulation. All material parameters and model dimensions are kept consistent with Figure \ref{fig:bench-ls-fric-setup2} and Table \ref{tab:bench-ls-fric-setup2}. The stress prediction of our XFEM model is compared to the Abaqus\textsuperscript{\textregistered} results in Figure \ref{fig:bench-ls-compare-vm2}. Furthermore, the results for the total contact force as a function of the applied displacements are shown in Figure \ref{fig:bench-ls-compare-cf2}. The relative difference between the XFEM and Abaqus\textsuperscript{\textregistered} force predictions, integrated over the loading process, is $6.5\times 10^{-4}$.

%===============================================================================

\subsection{Sensitivity Analysis Verification}
\begin{figure}[t]
\begin{minipage}[t]{0.42\textwidth}
\centering
\vspace{0pt}
  \includegraphics[width=1\linewidth]{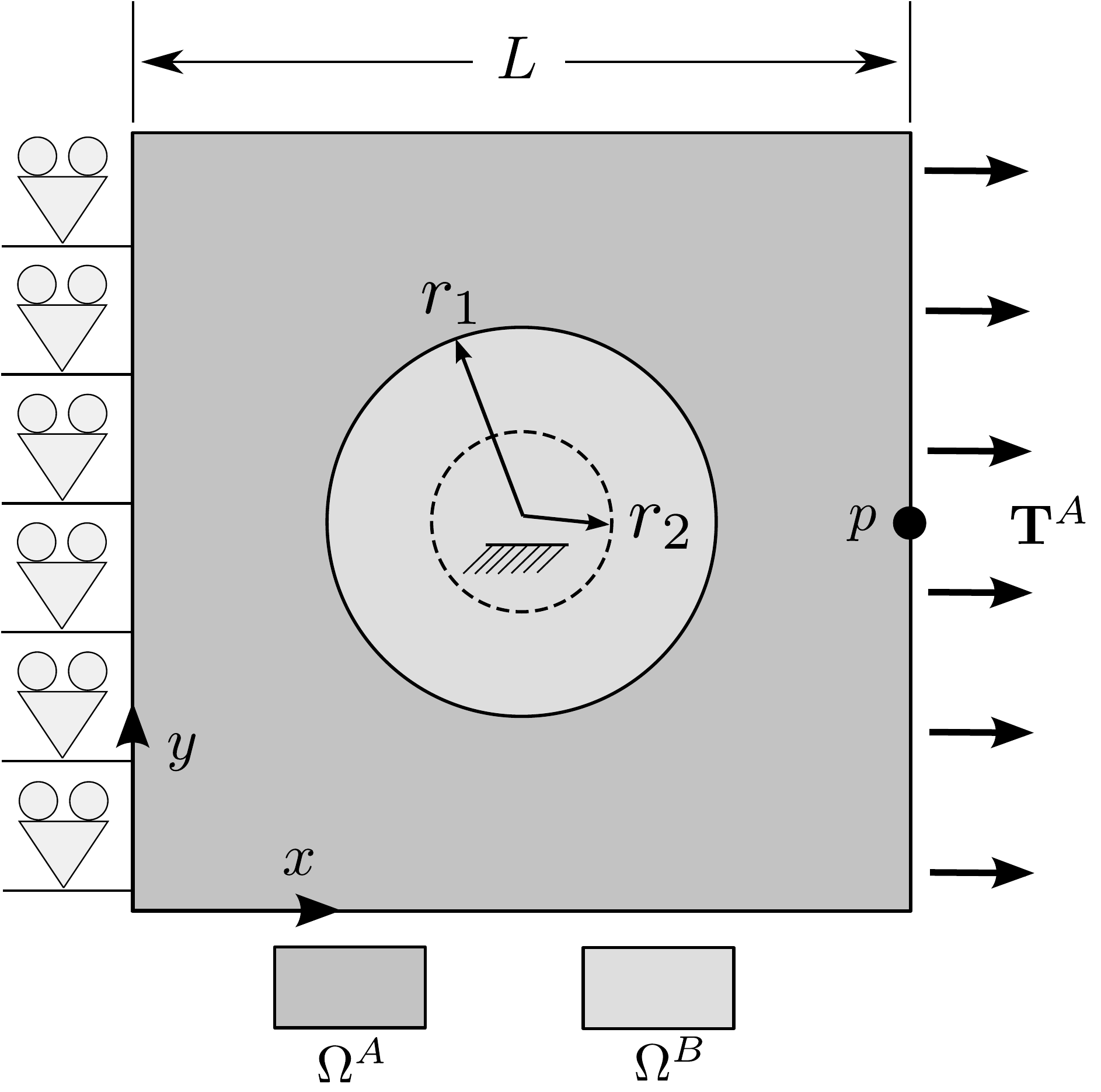}
  \caption{Sensitivity verification model setup.}
  \label{fig:sens_setup}
\end{minipage}\hfill
\begin{minipage}[t]{0.48\textwidth}
\centering
\vspace{0pt}
\captionof{table}{Sensitivity verification model parameters.}
\label{tab:sens_setup}
\begin{tabular}{ll}
\hline
Description & Parameter \\ \hline
domain length & $L = 1.0$ m \\
interface radius & $r_1 = 1.2$ m \\
fixed radius & $r_2 = 0.125$ m \\
Young's modulus & $E^A = 10$ MPa \\
Young's modulus & $E^B = 10$ MPa \\
Poisson's ratio & $\nu^A = 0.3$ \\
Poisson's ratio & $\nu^B = 0.3$ \\
applied load & $\mathbf{T}^A = (10,0)$  kPa/m \\
observation point & $p$ \\
\hline
\end{tabular}
\end{minipage}
\end{figure}
\begin{figure}[t]
\begin{minipage}[t]{0.48\textwidth}
\centering
\vspace{0pt}
  \includegraphics[width=.9\linewidth]{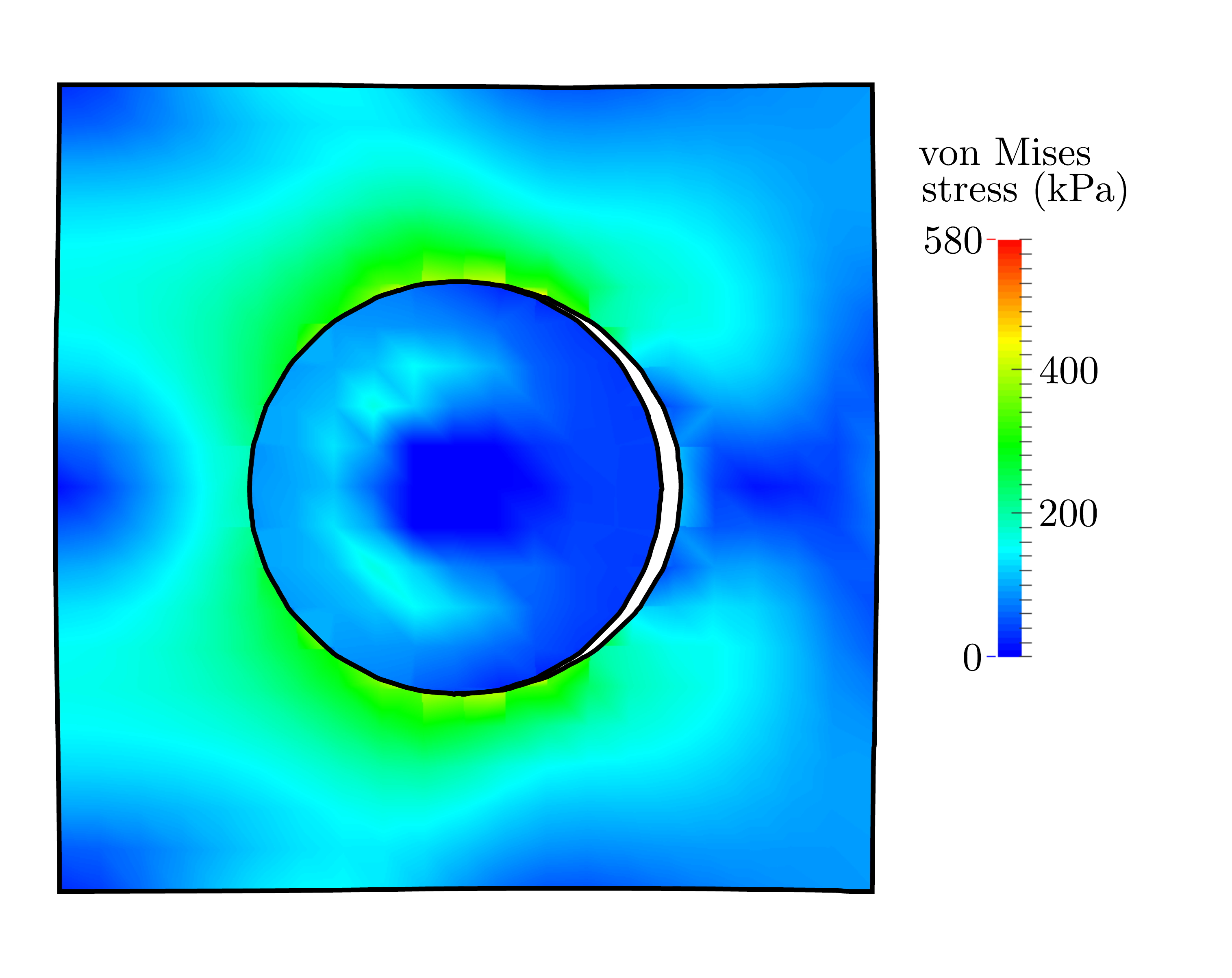}
  \caption{Von Mises stress distribution for $r_1 = 0.275$ m.}
  \label{fig:sens_nom_resp}
\end{minipage}\hfill
\begin{minipage}[t]{0.50\textwidth}
\vfill
\centering
\vspace{0pt}
\includegraphics[width=1\linewidth]{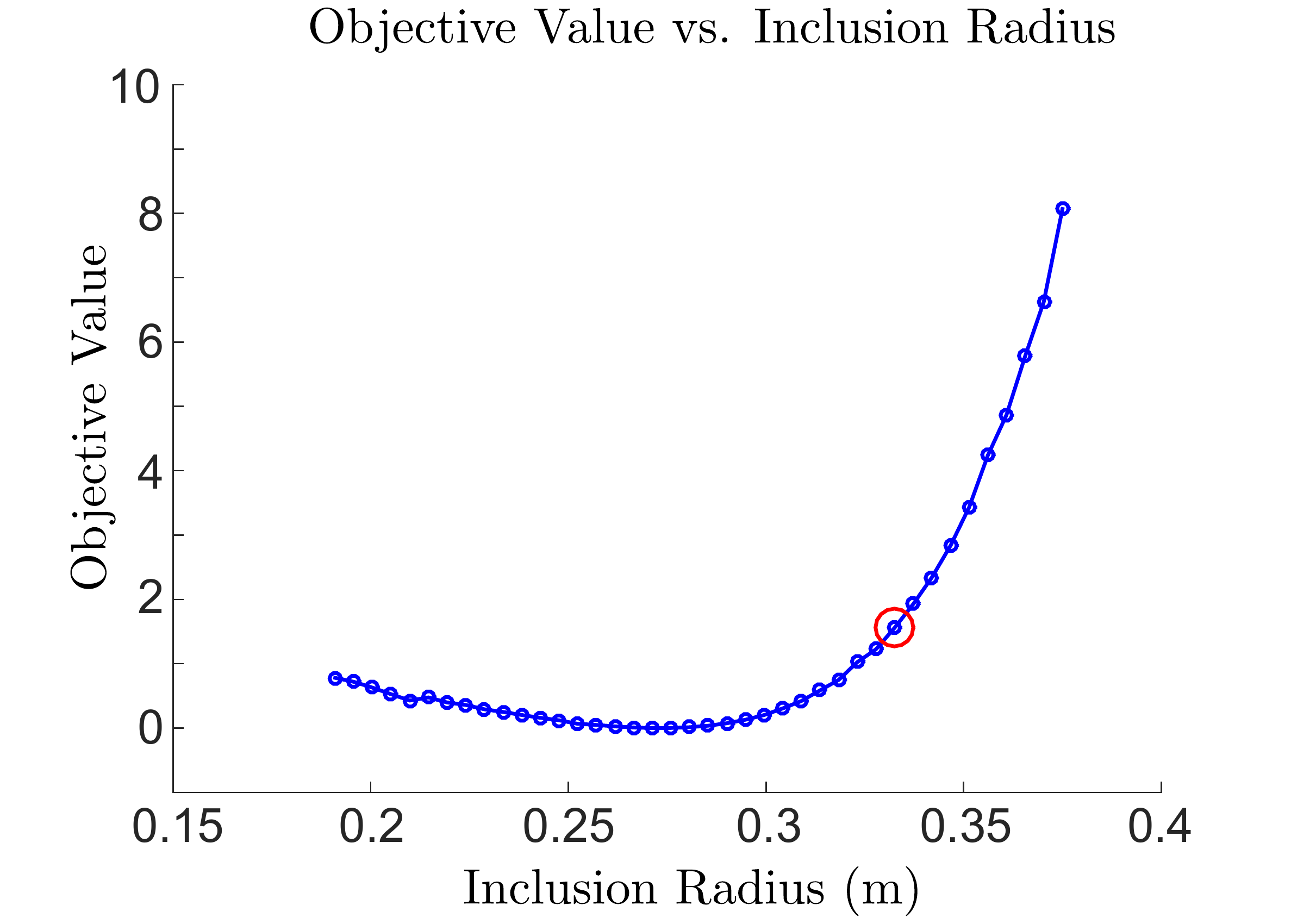}
\caption{Objective value as a function of material interface radius; sensitivity evaluated for range of radii highlighted in red.}
\label{fig:sens_obj_sweep}
\end{minipage}
\end{figure}
To verify the accuracy of the design sensitivities evaluated by the adjoint method for problems involving large sliding contact, we consider the optimization problem illustrated in Figure \ref{fig:sens_setup}. The square design domain of length $L$ is held in place by a circular inclusion of radius $r_1$, fixed within radius $r_2$. In addition, the left hand edge of the design domain is constrained in the $y$-direction. An external traction $\mathbf{T}^A$ is distributed along the right hand edge and applied in two load steps. Frictionless contact is modeled at the material interface. The converge criterion for solving the nonlinear systems in each load step by Newton's method is set to $1\times10^{-9}$. The XFEM model is discretized with 20$\times$20 elements. The model parameters are listed in Table \ref{tab:sens_setup}.

The objective function is defined as:
\begin{equation}
z = 1\times10^{4}\left(u_x\left(p\right) - 2.7731\times10^{-2}\right)^2.
\end{equation}
The design variable, $s$, defines radius of the circular inclusion, i.e.~$r_1=s$. The inclusion is described by the following LSF:
\begin{equation}
\phi = r_1 - \sqrt{\left(\mathbf{x} - 0.5\right)^2 + \left(\mathbf{y} - 0.5\right)^2}.
\end{equation}
To evaluate the behavior of the objective function with respect to the design variable, the inclusion radius is swept from $0.191\leq r_1 \leq 0.375$. Figure \ref{fig:sens_nom_resp} illustrates the mechanical response at $r_1 = 0.275$ m. The objective value over the interface radius is plotted in Figure \ref{fig:sens_obj_sweep}. The results show a rather smooth dependency of the objective on the interface radius.
\begin{figure}[t]
\begin{center}
\includegraphics[width=.65\linewidth]{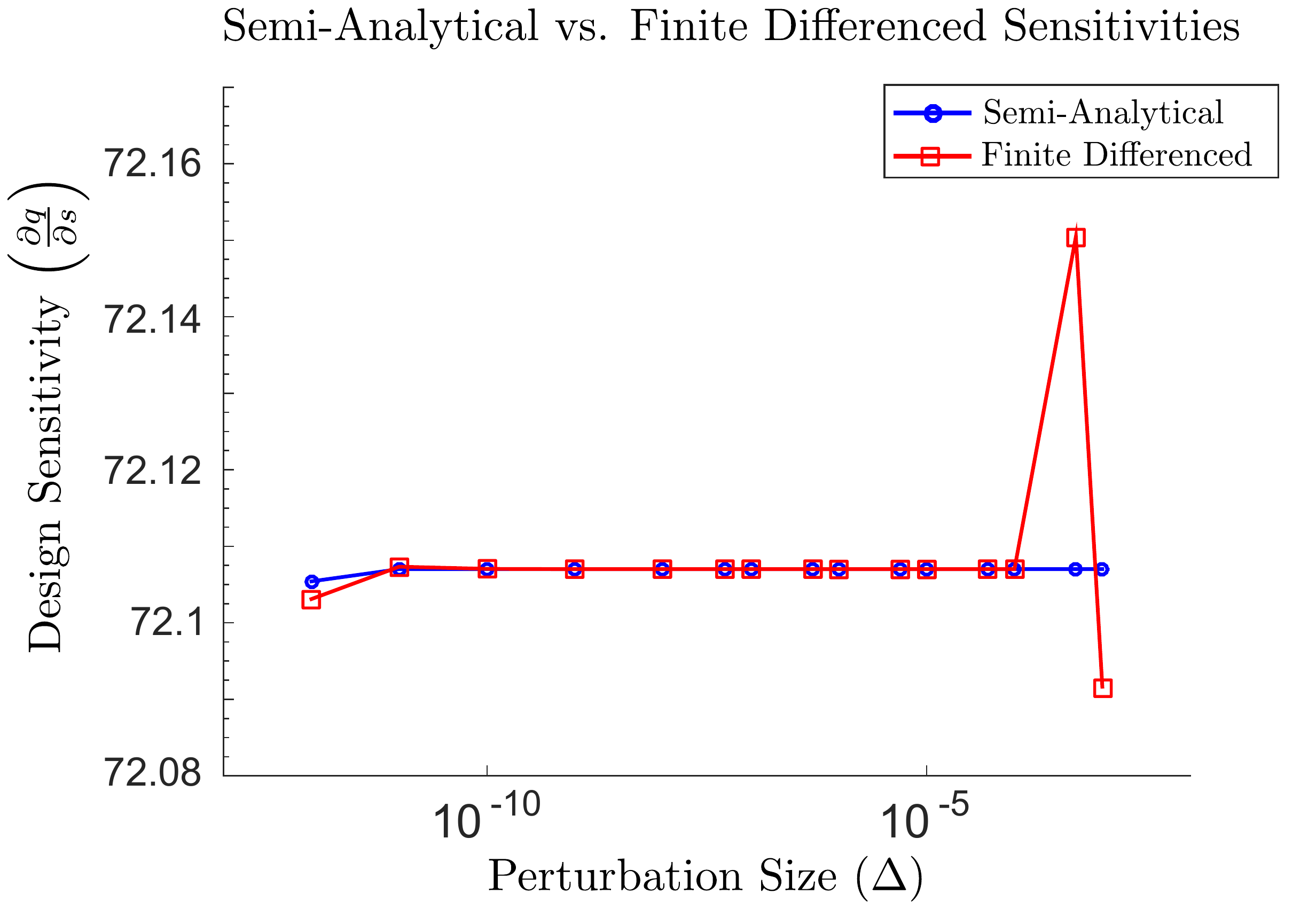}
\caption{Sensitivity of objective value with respect to inclusion radius evaluated by semi-analytical adjoint method and finite differencing over a range of perturbation sizes.}
\label{fig:fd_vs_sa}
\end{center}
\end{figure}

To verify our adjoint sensitivity analysis method, the semi-analytical evaluation of the sensitivities is compared to the results of a central finite differencing scheme.
%
%\begin{equation}
%\frac{\mathrm{d} z}{\mathrm{d} s} = \frac{\partial z}{\partial s} + \frac{\partial z}{\partial \mathbf{u}}{\frac{\partial \mathbf{R}}{\partial \mathbf{u}}}^{-1}\left[-\frac{\Delta \mathbf{R}}{\Delta s}\right]
%\end{equation}
%
%where $\Delta$ denotes the perturbation of a quantity.
%
Figure \ref{fig:fd_vs_sa} plots the semi-analytical and finite differenced sensitivities as a function of the perturbation size. At extremely small perturbation sizes, machine precision round off errors affect the accuracy of the design sensitivities. For this particular problem, perturbation sizes larger than $\Delta s > 10^{-4}$ yield linearization errors when using finite differencing. For perturbations between $10^{-8}$ and $10^{-5}$ there is good agreement between both semi-analytical and finite differenced values. The relative error between both methods at a perturbation size of $10^{-8}$ is $4.62\times10^{-10}$. However, in addition to being computationally less expensive, the adjoint method is less sensitive to the perturbation size than finite differenced design sensitivities.

%===============================================================================

\subsection{Material Anchor Design Problem}\label{sec:anchor}

To compare designs optimized with finite strain theory to designs optimized with infinitesimal strain theory, we study a material anchor design problem.  A structural anchor is embedded within a host material of the same properties, with frictionless contact at the interface to afford resistance to separation. The objective is to find an optimized anchor geometry such that the resistance is maximized. This problem was originally studied under a small strain assumption by \cite{LM:15}.
\begin{figure}[t]
\begin{minipage}[t]{0.47\textwidth}
\centering
\vspace{0pt}
\includegraphics[width=1\linewidth]{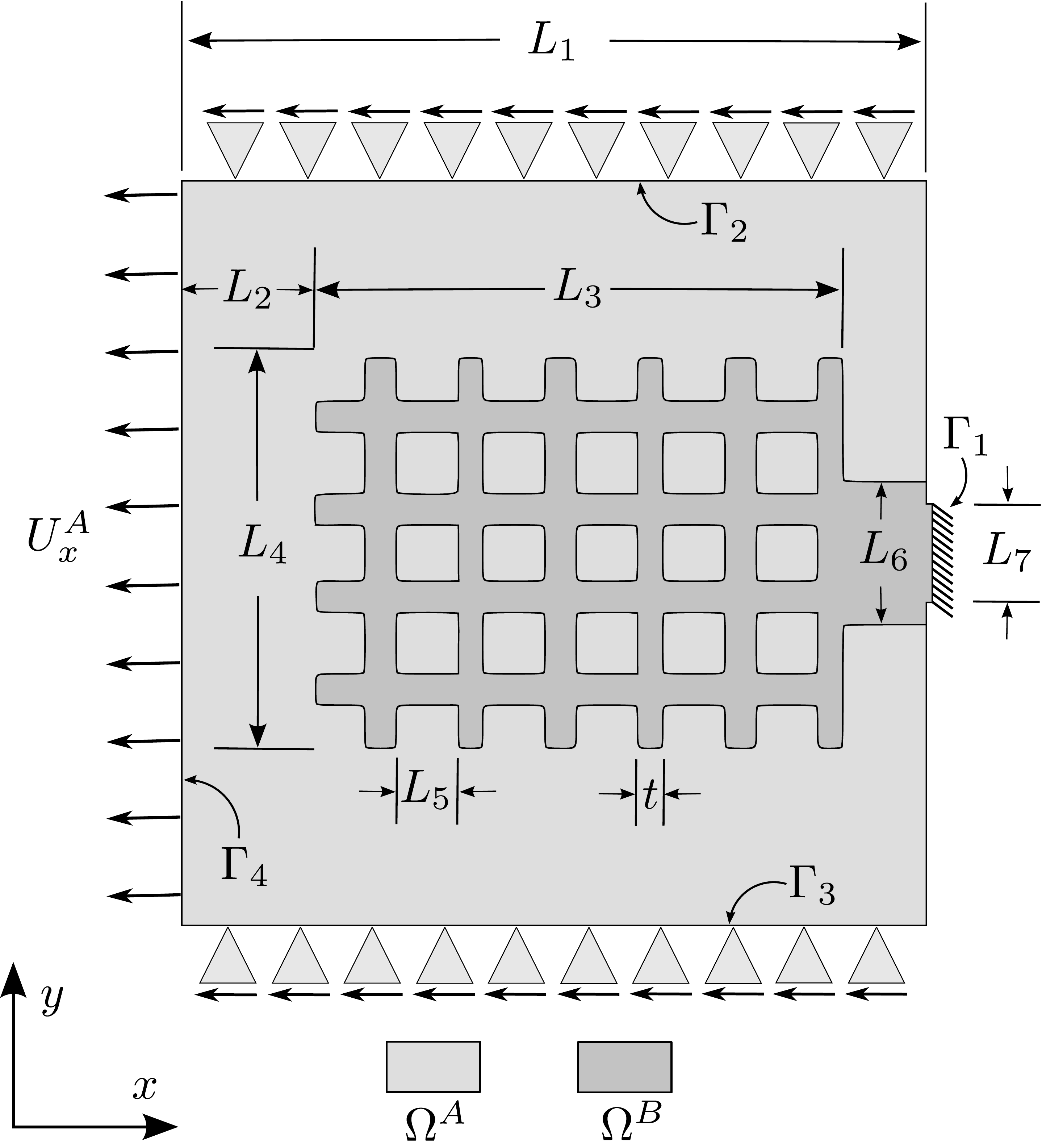}
\caption{Material anchor initial configuration.}
\label{fig:ex1_setup}
\end{minipage}\hfill
\begin{minipage}[t]{0.51\textwidth}
\centering
\vspace{0pt}
\captionof{table}{Nominal material anchor model parameters.}
\label{tab:ex1_setup}
\begin{tabular}{ll}
\hline
Description & Parameter \\ \hline
domain length & $L_1  = 1.0$ m \\
host depth & $L_2     = 0.175$ m \\
lattice length & $L_3 = 0.715$ m \\
lattice width & $L_4  = 0.523$ m \\
cuboid length & $L_5  = 0.0833$ m \\
anchor base width & $L_6   = 0.2$ m \\
fixed support width & $L_7 = 0.133$ m \\
lattice thickness & $t = 0.038$ m \\
Young's modulus & $E^A = 10$ MPa \\
Young's modulus & $E^B = 10$ MPa \\
Poisson's ratio & $\nu^A = 0.3$ \\
Poisson's ratio & $\nu^B = 0.3$ \\
applied displacement & $U_x^A = 0.01$ m \\
response weight &  $c_u = 75$ \\
penalty weight & $c_p   = 25$ \\
volume ratio & $c_v     = 0.5$ \\
opt. upper bounds & $s_{max} = 8.33\times 10^{-3}$ \\
opt. lower bounds & $s_{min} = -8.33\times 10^{-3}$ \\
rel. step size & $\Delta s = 8\times10^{-3}$\\
smoothing radius & $r_f = 0.0375$ m \\
\hline
\end{tabular}
\end{minipage}
\end{figure}

The initial material distribution and boundary conditions are illustrated in Figure \ref{fig:ex1_setup}, while model parameters are listed in Table \ref{tab:ex1_setup}. The volume occupied by the anchor material, $\Omega^B$, is fixed at $\Gamma_1$, while a prescribed displacement, $U_x^A$, is applied to the volume occupied by the host material, $\Omega^A$, along $\Gamma_{2-4}$. Displacements are constrained to zero in the $y$ direction along $\Gamma_{2-4}$. To prevent the anchor material from directly connecting boundary $\Gamma_1$ to boundaries $\Gamma_{2-4}$, they are excluded from the design domain.

For this example, we wish to determine the optimal geometry such that the force at $\Gamma_1$ is maximized. The mechanical response contribution to the objective function in Equation~\eqref{eq:opt} is defined as:
\begin{equation}
z = 100 - \int_{\Gamma_1}\sigma_{xx}\mathrm{d}\Gamma,
\end{equation}
where $\sigma_{xx}$ denotes the normal Cauchy stress in the $x$ direction. To regularize the problem, a perimeter penalty of $c_p = 0.25$ is applied. To prevent the anchor material from occupying the majority of the design space, a volume constraint of 50 \% is applied to the anchor material, i.e. $c_v = 0.5$. Due to the symmetric nature of the problem, only one half of the domain is analyzed with 120$\times$60 elements. The linear level set filter~\eqref{eq:smooth_rad} for this problem is set to $r_f = 4.5h$, where $h$ is the element side length.

\subsubsection{Nominal Design}
\label{sec:nomdes}
\begin{figure}[t]
\begin{center}
\includegraphics[width=.85\linewidth]{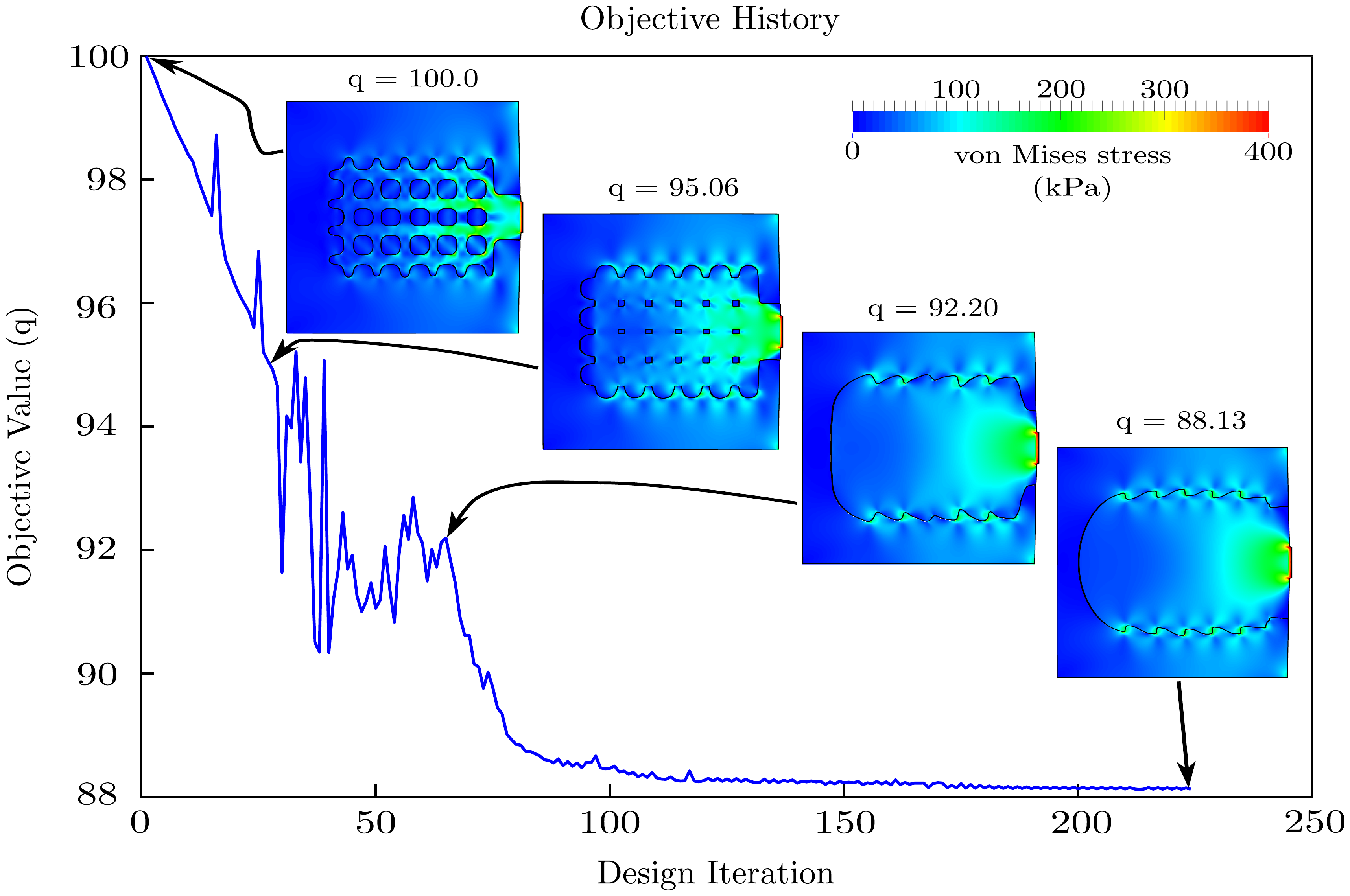}
\caption{Material anchor design objective history with snapshots of specific iterations. Inset depicts the force-displacement curve for specific iterations.}
\label{fig:ex1_final}
\end{center}
\end{figure}

We first present the results of the material anchor nominal design, using the design parameters listed in Table \ref{tab:ex1_setup}. The magnitude of applied load for this example, $U_x^A = 0.01$ m, was specifically chosen to keep the experienced strain well within the limitations of small strain theory. This allows comparing the results of the proposed optimization method for large strain contact with the results of \cite{LM:15} where infinitesimal strains, a linear elastic response, and negligible sliding between surfaces were assumed. Here, the displacement, $U_x^A$, is applied in two load steps.

Figure \ref{fig:ex1_final} illustrates the objective value history during optimization, supported by snapshots of the mechanical response for specific design iterations. The anchor material quickly merges to a uniform body, producing ridges or spines along the outer surface to afford resistance to separation. In the early stages of convergence, topological changes result in an abrupt change in the measured objective value. Once the topology remains unchanged, the optimization process converges smoothly. The optimized geometry closely resembles the small strain theory analog presented by \cite{LM:15}.

\subsubsection{Load Case Study}

Even for problems which maintain small relative motion between surfaces in contact, infinitesimal strain theory may still result in mechanical response evaluation inaccuracies. Infinitesimal strain theory assumes that the surface orientation in the current configuration differs insignificantly from that of the undeformed configuration. This assumption may lead to errors in contact pressure estimation and affect the optimized geometry. To illustrate this issue, the applied displacement is increased, and the material anchor problem is optimized for both large and small strain theory. For this load case study, we adopt a continuation approach, where the second load case uses the previous optimized geometry as an initial configuration.
\begin{figure}[t]
\begin{center}
\includegraphics[width=.7\linewidth]{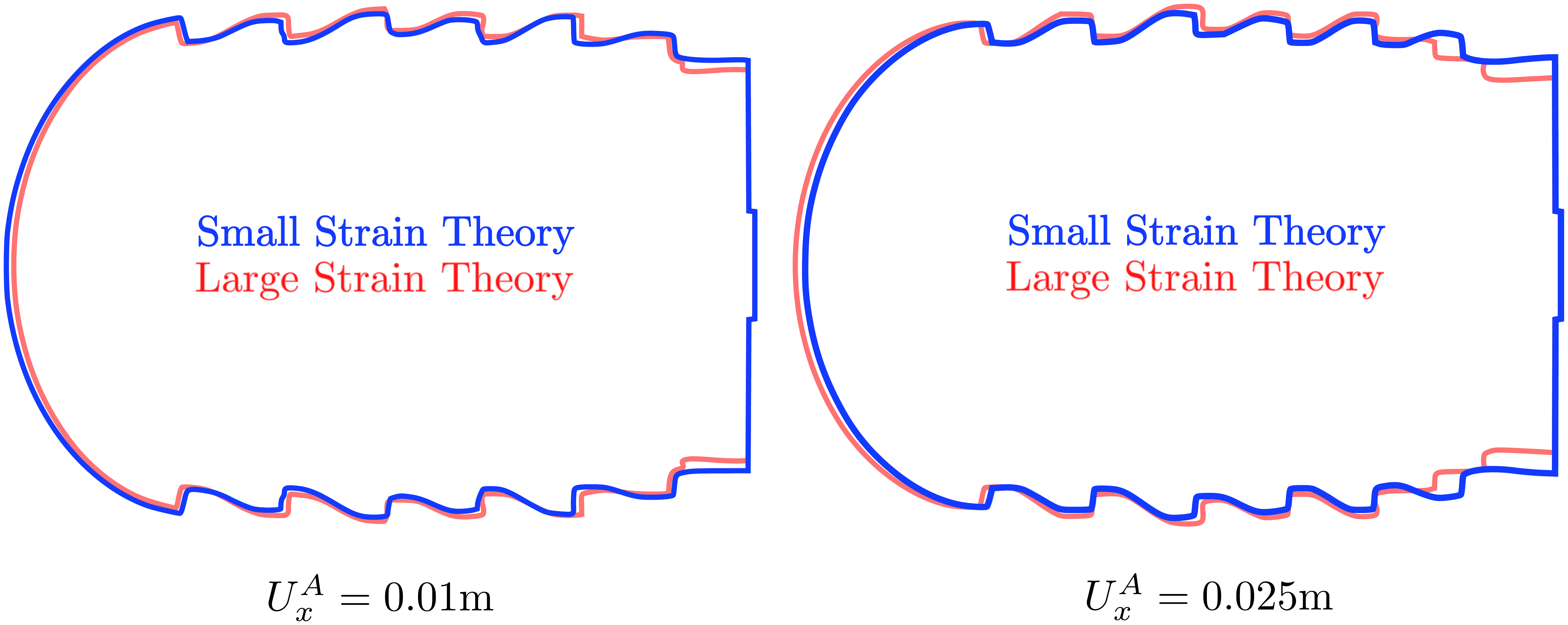}
\caption{Comparison of optimized geometry for small and large strain theory for various applied displacements.}
\label{fig:ex1_compare}
\end{center}
\end{figure}
\begin{figure}[b]
\begin{center}
\includegraphics[width=.55\linewidth]{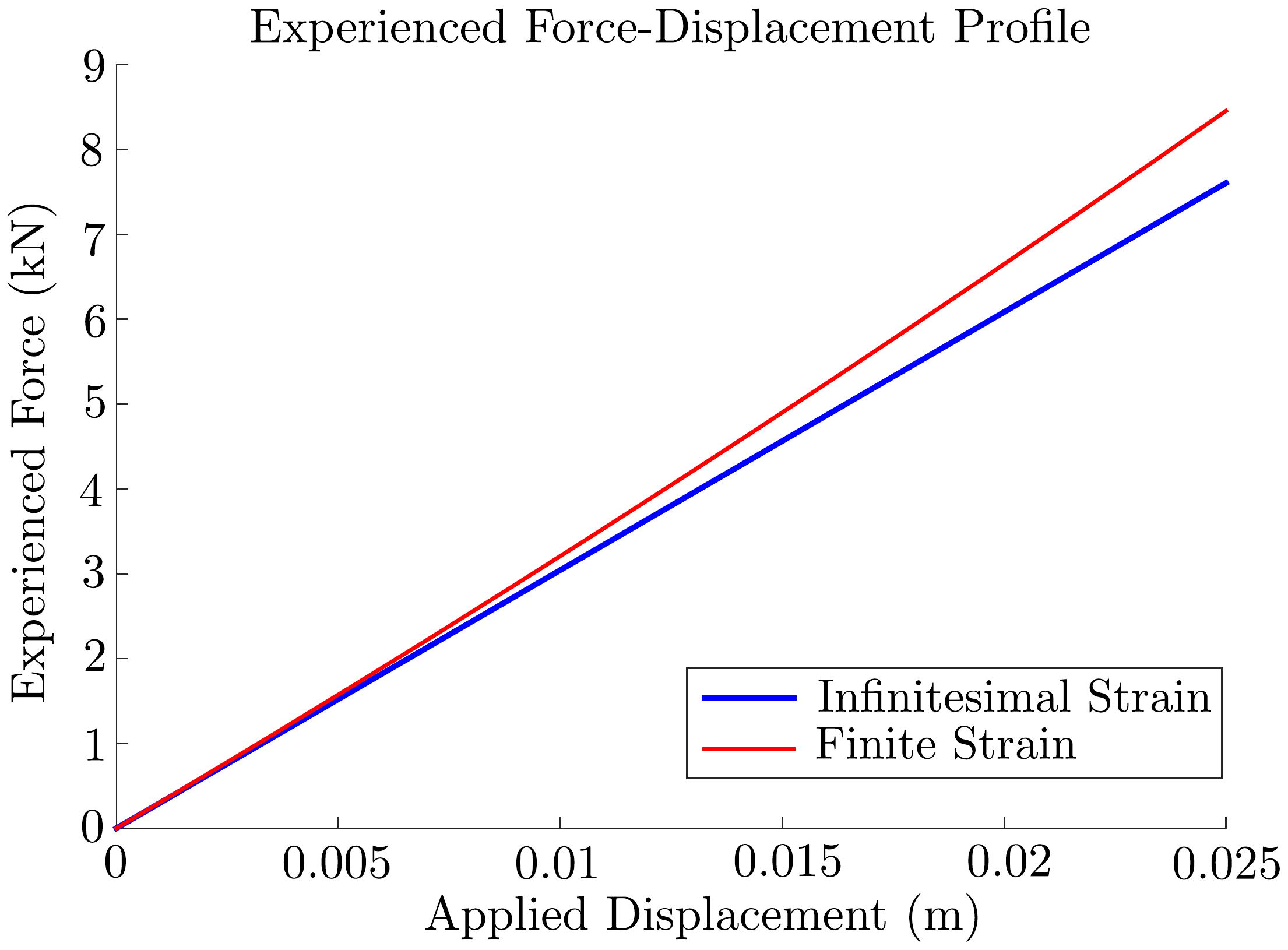}
\caption{Comparison of resistance to separation for incremental loading to the maximum value $U^A_x = 0.025$ m.}
\label{fig:ex1_fd_compare}
\end{center}
\end{figure}

Figure \ref{fig:ex1_compare} compares small and large strain theory optimized geometries for two different magnitudes of applied loads. At the smallest load case, $U^A_x = 0.01$ m, the optimized material anchor profile produced by small strain theory closely resembles that resolved by large strain theory. However, at the higher load case, $U^A_x = 0.025$ m, the discrepancy between the optimized geometries becomes more noticeable, although the conceptual designs differ insignificantly. The close resemblance in geometry can be attributed to similar physical responses. To illustrate these similarities, the optimal geometries for both infinitesimal strain and finite strain theory are incrementally loaded in $15$ steps for the larger load case of $U^A_x = 0.025$ m. Figure \ref{fig:ex1_fd_compare} presents the experienced force profile as a function of applied displacement. At this load level, the large strain contact model displays a rather linear response, similar to the infinitesimal strain contact model.
\begin{table}
\centering
\def\arraystretch{1.15}
\begin{tabular}{ccc}
\hline
 Optimized geometry for: &  \multicolumn{2}{c}{Holding force for load case:}  \\
               &  $U^A_x = 0.01$ m    & $U^A_x = 0.025$ m      \\
\hline
Finite strain    & $3.155$  kN         & $8.456$ kN  \\
Infinitesimal strain        & $3.1516$ kN        & $8.309$ kN   \\
\hline \vspace{0.03cm}
\end{tabular}
\caption{Holding force of optimized geometries using finite strain theory.}
\label{table:ex1_cross_compare}
\end{table}

To cross examine the performance of geometries optimized with infinitesimal strain theory and finite strain theory, the optimized geometries are analyzed with finite strain theory. Table \ref{table:ex1_cross_compare} compares the holding force of the geometries provided in Figure \ref{fig:ex1_compare} when analyzed strictly with finite strain theory. For the loading cases $U^A_x = 0.01$ m and $U^A_x = 0.025$ m, the optimized geometry produced from small strain theory results in a decrease in holding force of $0.11\%$ and $1.73\%$ respectively, when compared to the optimized geometry produced from large strain theory. While the improved performance of optimal geometry from finite strain theory is relatively small, this example demonstrates the limitations of infinitesimal strain theory for these types of optimization problems. To explore optimization problems in which the physical behavior cannot be predicted with any acceptable accuracy using infinitesimal strain theory, the following examples study design problems with highly nonlinear response behavior.

%===============================================================================

\subsection{Snap-Fit Design Problem}\label{sec:snapfit}
\begin{figure}[b]
\begin{minipage}[t]{0.47\textwidth}
\centering
\vspace{0pt}
\includegraphics[width=1\linewidth]{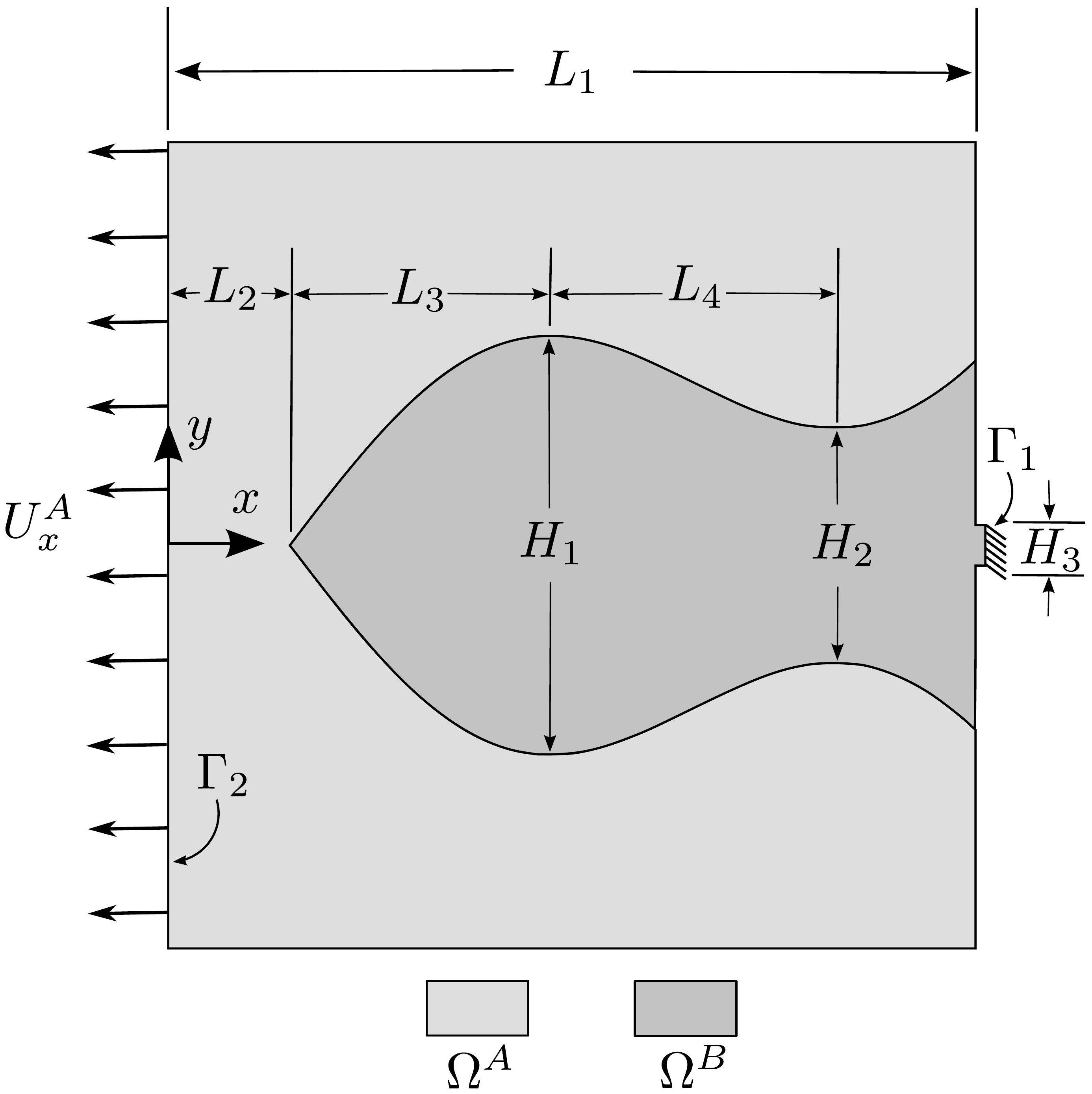}
\caption{Snap-fit design initial configuration.}
\label{fig:ex2_setup}
\end{minipage}\hfill
\begin{minipage}[t]{0.51\textwidth}
\centering
\vspace{0pt}
\captionof{table}{Snap-fit design model parameters.}
\label{tab:ex2_setup}
\begin{tabular}{ll}
\hline
Description & Parameter \\ \hline
domain length & $L_1 = 1.0$ m \\
host depth & $L_2 = 0.151$ m \\
peak width location & $L_3 = 0.4$ m \\
base width location & $L_4 = 0.9$ m \\
peak height & $H_1 = 0.25$ m \\
base height & $H_2 = 0.16$ m \\
fixed support height & $H_3 = 0.133$ m \\
Young's modulus & $E^A = 10$ MPa \\
Young's modulus & $E^B = 10$ MPa \\
Poisson's ratio & $\nu^A = 0.3$ \\
Poisson's ratio & $\nu^B = 0.3$ \\
applied load, at $t = 1$ & $U_x^A = 0.5$ m \\
response weight &  $c_u = 100.0$ \\
penalty weight & $c_p = 0.0$ \\
volume ratio & $c_v = 1.0$ \\
opt. upper bounds & $s_{max} = 0.0125$ \\
opt. lower bounds & $s_{min} = -0.0125$ \\
rel. step size & $\Delta s = 8\times10^{-3}$\\
smoothing radius & $r_f = 0.0375$ m \\
\hline
\end{tabular}
\end{minipage}
\end{figure}

Snap-fits remain one of the fastest and cost effective methods of assembly. This simplistic fastener relies on two interlocking components, which if designed properly can be assembled and disassembled numerous times without damaging the components. For applications demanding a high level of precision, the force required to induce separation can be pivotal. Snap-fit designs exhibit a highly nonlinear mechanical response, and during the process of separation the mechanical model can experience moments of neutral and unstable equilibrium. This poses interesting challenges for both the mechanical response prediction, and subsequently the optimization of snap-fit designs.

Here we pose the snap-fit design problem as follows: we wish to find the optimal geometry of a snap-fit mechanism to match a desired load-displacement profile. This problem formulation is explored for two scenarios: the first example is a two phase design in which geometry control is provided by discretized level set nodal variables; the second example is a three-phase design in which geometry control is provided by geometric primitive variables.

%===============================================================================

\subsubsection{Two-Phase Example}
\begin{figure}[t]
\begin{center}
\includegraphics[width=1\linewidth]{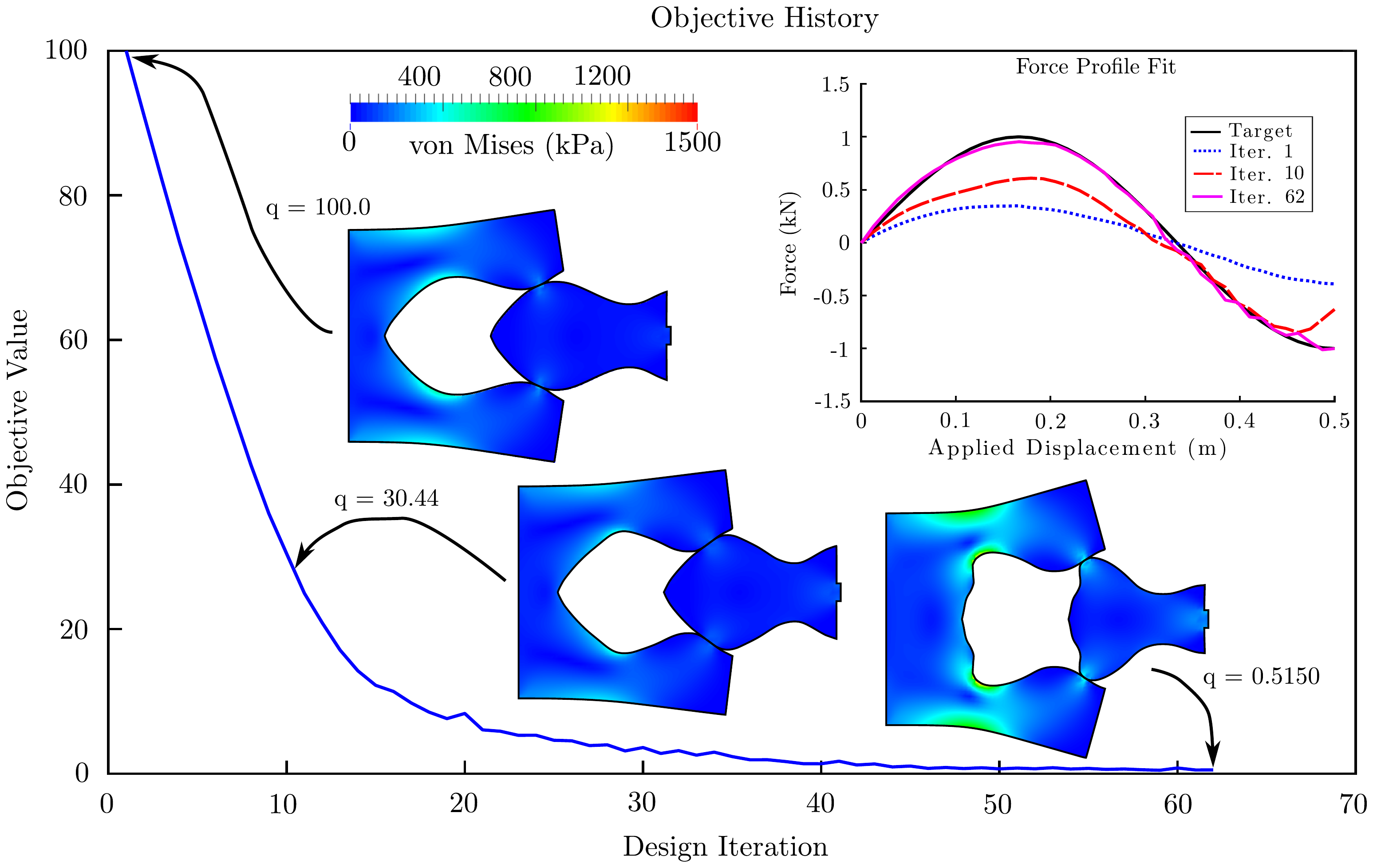}
\caption{Snap-fit design objective history with snapshots of specific iterations. Inset depicts the force-displacement curve for specific iterations.}
\label{fig:ex2_final}
\end{center}
\end{figure}

The initial configuration is illustrated in Figure \ref{fig:ex2_setup}, and dimensions and material parameters are listed in Table \ref{tab:ex2_setup}. The snap-fit tab, represented by Phase B, is fixed at the boundary $\Gamma_1$, whereas the snap-fit container, represented by phase A, is subjected to a prescribed displacement along boundary $\Gamma_2$, which is applied in $45$ equal load increments.

The objective function is defined as follows:
\begin{equation}\label{eq:ex2_obj}
z = \int \left( \int_{\Gamma_1} \sigma_{xx}\mathrm{d}\Gamma  - f_t \right)^2\mathrm{d}t,
\end{equation}
where $f_t$ is a target force profile, and the pseudo-time $t$ represents the incremental loading process. For this particular example, the desired force profile is defined as:
\begin{equation}
f_t = \sin\left(\frac{3\pi}{2} t \right) \ \mathrm{kN}
\qquad 0\leq t \leq 1.
\end{equation}
The desired force profile describes a sinusoidal profile with a  peak value of 1 kN. This particular force profile was chosen to highlight a design exhibiting a high level of physical response non-linearity.

For geometry control, the nodal level set values are defined in terms of the optimization variables using the linear filter \eqref{eq:smooth_rad}. The smoothing radius is set to $r_f=3h$, where $h$ is the element side length. No perimeter penalty measures or volume constraints are used for this example. The smoothness of the target force-profile causes localized geometry irregularities to be non-beneficial to the design functionality. Due to the symmetric nature of the design, only half of the domain is analyzed with 48$\times$24 elements.

Figure \ref{fig:ex2_final} shows the convergence profile of the optimization problem, supported by snapshots of the mechanical response at the final time step for select design iterations. The inset of Figure \ref{fig:ex2_final} shows the experienced force-displacement profile, compared to the desired profile for specific design iterations. The stem of the tab increases in concavity, increasing the experienced force at the base. The peak width of the tab is increased, whereas the pointed head of the tab flattens out. The geometry evolution observed increases the peak force experienced, and provides a close fit to the desired force-displacement profile. The non-smooth nature of the final design force-displacement curve can be attributed to the piecewise linear interface representation in the XFEM model.

%===============================================================================

\subsubsection{Three-Phase Example}
\label{sec:3psfit}
\begin{figure}[t]
\begin{minipage}[t]{0.52\textwidth}
\centering
\vspace{0pt}
\includegraphics[width=1\linewidth]{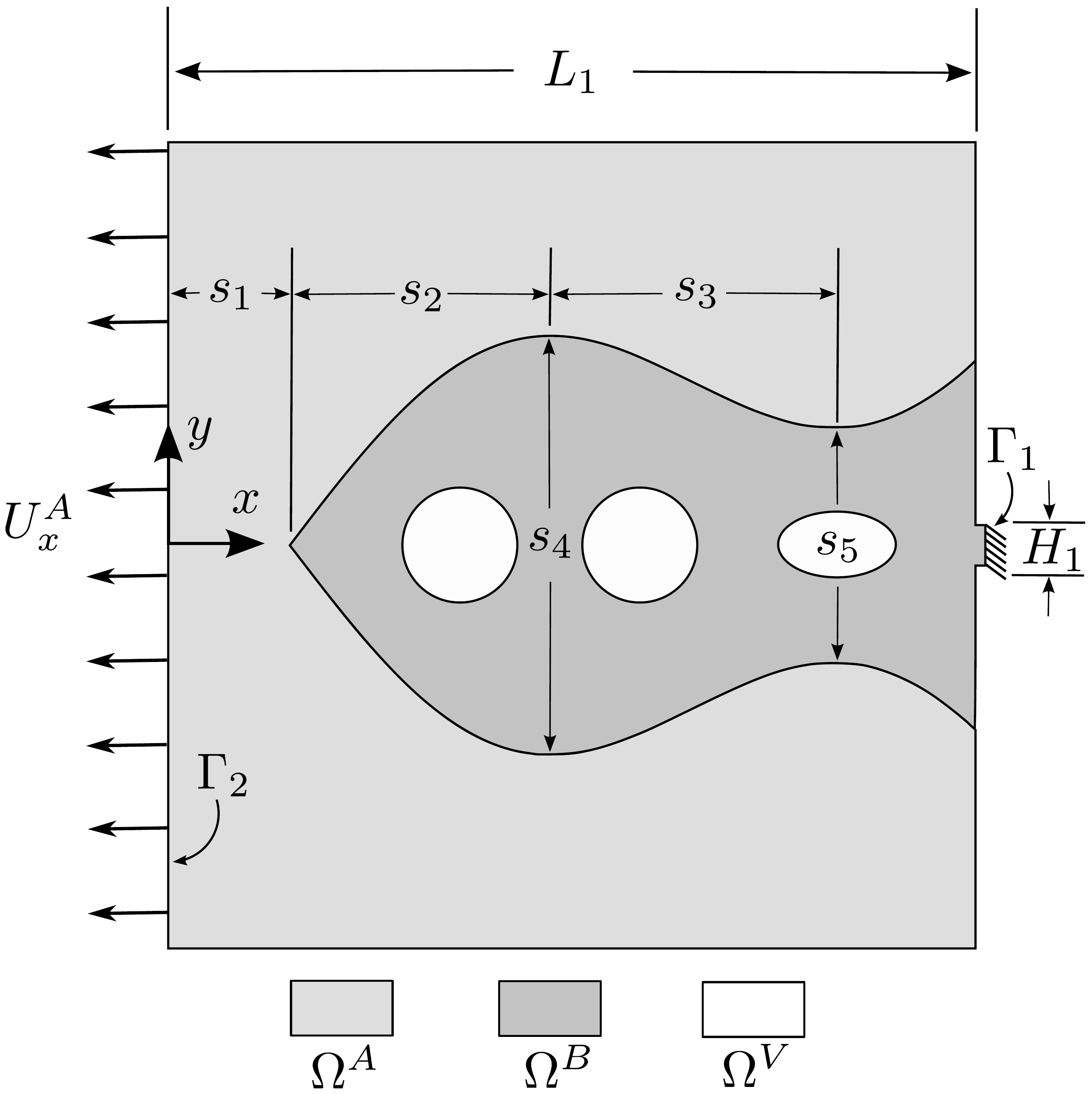}
\caption{Snap-fit design initial configuration.}
\label{fig:ex2_3p_setup}
\end{minipage}\hfill
\begin{minipage}[t]{0.47\textwidth}
\centering
\vspace{0pt}
\captionof{table}{Snap-fit design model parameters.}
\label{tab:ex2_3p_setup}
\begin{tabular}{ll}
\hline
Description & Parameter \\ \hline
domain length & $L_1 = 1.0$ m \\
fixed support height & $H_1 = 0.133$ m \\
host depth & $s_1 = 0.151$ m \\
peak width location & $s_2 = 0.4$ m \\
base width location & $s_3 = 0.9$ m \\
peak height & $s_4 = 0.25$ m \\
base height & $s_5 = 0.16$ m \\
Young's modulus & $E^A = 10$ MPa \\
Young's modulus & $E^B = 10$ MPa \\
Poisson's ratio & $\nu^A = 0.3$ \\
Poisson's ratio & $\nu^B = 0.3$ \\
applied displacement & $U_x^A = 0.6$ m \\
response weight &  $c_u = 99.9$ \\
penalty weight & $c_p = 0.1$ \\
volume ratio & $c_v = 0.15$ \\
rel. step size & $\Delta s = 4\times10^{-4}$\\
\hline
\end{tabular}
\end{minipage}
\end{figure}

The two-phase snap-fit design demonstrated the optimization of a problem experiencing large sliding motion, but with rather small elastic deformations. The three-phase analog explores an optimization problem exhibiting large sliding contact in the presence of large elastic deformation, by introducing void regions within the tab material. The three-phase snap-fit design problem describes the design domain using a combination of geometric primitives, wherein the variables associated with the geometric primitives are defined as optimization variables. While the design freedom is restricted to the set of shapes that can be produced by the particular geometric primitives, the three-phase example explores a different avenue of complexity by introducing void regions within the tab material, i.e.~phase B.

The initial configuration is illustrated in Figure \ref{fig:ex2_3p_setup}, and dimensions and material parameters are listed in Table \ref{tab:ex2_3p_setup}. Similar to the two-phase snap-fit example, the tab material is fixed at the boundary $\Gamma_1$ and a prescribed displacement is applied at boundary $\Gamma_2$ in $45$ load steps. Provided the same objective of matching a target force displacement profile measured at $\Gamma_1$, a new target force displacement curve is defined as:
\begin{equation}
f_t = 0.5 \sin\left(\frac{3\pi}{2} t \right) \ \mathrm{kN}
\qquad 0\leq t \leq 1.
\end{equation}
Provided the increased applied displacement value of $U_x^A = 0.6$ m, the target force displacement curve for this example exhibits a delayed peak force value as compared to the two-phase example. This is done to promote larger deformations prior to exceeding the point of neutral equilibrium.
\begin{table}[t]
\centering
\begin{tabular}{lllll}
\hline
Description & Variable & Initial Value & Upper Bound & Lower Bound \\ \hline
host depth           & $s_1$     & $0.151$ m & $0.3$ m  & $0.1$ m \\
peak width location  & $s_2$     & $0.4$ m   & $0.7$ m  & $0.25$ m  \\
base width location  & $s_3$     & $0.9$ m   & $1.0$ m  & $0.8$ m \\
peak height          & $s_4$     & $0.2451$ m  & $0.4$ m  & $0.15$ m \\
base height          & $s_5$     & $0.16$ m  & $0.25$ m & $0.04$ m  \\
$\phi_{c,1}$ x center      & $x_{c,1}$ & $0.4$ m   & $0.85$ m & $0.25$ m  \\
$\phi_{c,2}$ x center      & $x_{c,2}$ & $0.6$ m   & $0.85$ m & $0.25$ m \\
$\phi_{c,3}$ x center      & $x_{c,3}$ & $0.8$ m   & $0.85$ m & $0.25$ m \\
$\phi_{c,1-3}$ x radius      & $r_{x,1-3}$ & $0.07$ m  & $0.3$ m  & $0.02$ m \\
$\phi_{c,1-2}$ y radius      & $r_{y,1-2}$ & $0.07$ m  & $0.2$ m  & $0.02$ m\\
$\phi_{c,3}$ y radius      & $r_{y,3}$ & $0.035$ m & $0.06$ m & $0.02$ m\\
\hline
\vspace{0.1cm}
\end{tabular}
\caption{Initial value, upper and lower bounds for three-phase snap-fit design problem.}\label{tab:ex2_3p_vars}
\end{table}

The design geometry is defined by two LSFs; see Section \ref{sec:3phasegeom}. The first LSF governs the material interface between subdomains $\Omega^A$ and $\Omega^B$ and is defined the optimization variables $s_1 - s_5$:
\begin{equation}
\phi^1 = - |Y| + a \ \sin\left(\theta\right) + p \ \tilde{X},
\end{equation}
where
\begin{equation}
a = \frac{s_4}{2} - \frac{p \ s_2}{s_2+s_3}, \qquad p = \frac{\left(s_4 + s_5\right)\left(s_2 + s_3\right)}{4 s_2 + 2 s_3},  \qquad \tilde{X} = \frac{X - s_1}{s_2+s_3},
\end{equation}
and the auxiliary coordinates,
\begin{equation}
\theta = \tilde{a}X^2 + \tilde{b}X + \tilde{c},
\end{equation}
are defined by the following scalar values
\begin{equation}
\tilde{a} = \frac{\pi\left(2 s_2 - s_3\right)}{2 s_2 s_3 \left(s_2 + s_3 \right)},
\end{equation}
\begin{equation}
\tilde{b} = \frac{\pi \left(\left(s_1 + s_2 + s_3\right)^2 - 3\left(s_1 + s_2\right)^2 + 2 s_1^2\right)}{s_2 s_3 \left(2 s_2 + 2 s_3\right)},
\end{equation}
\begin{equation}
\tilde{c} = -\frac{\pi \ s_1 \left(-2 s_2^2 + 2 s_2 s_3 - 2 s_1 s_2 + s_3^2 + s_1 s_3\right)}{s_2 s_3 \left(2 s_2 + 2 s_3\right)},
\end{equation}
The geometric primitives chosen for the LSF $\phi^1$ conveniently allow control over the length, peak width, and narrow width of the tab material outer geometry, while maintaining a smooth curvature of the profile.

The second LSF distinguishes $\Omega^V$ from $\Omega^B$. The geometric primitives chosen for the LSF $\phi^2$ define ellipsoidal void regions, shown in Figure \ref{fig:ex2_3p_setup}, which can move along the x-axis and grow or shrink in size. The LSF $\phi^2$ is defined as a combination of elliptical conical LSFs,
\begin{equation}
\phi_{c,j} = \left(\left(\frac{X - x_{c,j}}{r_{x,j}}\right)^2 + \left(\frac{Y - y_{c,j}}{r_{y,j}}\right)^2\right)^2 - 1,
\end{equation}
where the $jth$-elliptical conical field $\phi_{c,j}$ is defined by variables $x_{c,j}$ and $y_{c,j}$ which control the central location of the ellipse, and parameters $r_{x,j}$ and $r_{y,j}$ which control the semi-axis radii of the ellipse. Using three elliptical conical LSFs, the LSF $\phi^2$ is defined as:
\begin{equation}
\phi^2 = \min\left(\phi_{c,1},\phi_{c,2}, \phi_{c,3}\right).
\end{equation}
The y-location of the conical fields are restricted as $y_{c,j} = 0$, and are excluded from geometry control.

The initial values, upper limits, and lower limits for all design variables for this problem are presented in Table \ref{tab:ex2_3p_vars}. A volume constraint of $c_v = 0.15$ is enforced to reduce the material of $\Omega^B$, and a perimeter penalty weight of $c_p = 0.1$ used to regularize the void material interface. Due to the symmetric nature of the design, only half of the domain is analyzed with 80$\times$40 elements. No level set smoothing filters are used in this example.
\begin{figure}[t]
\begin{center}
\includegraphics[width=1\linewidth]{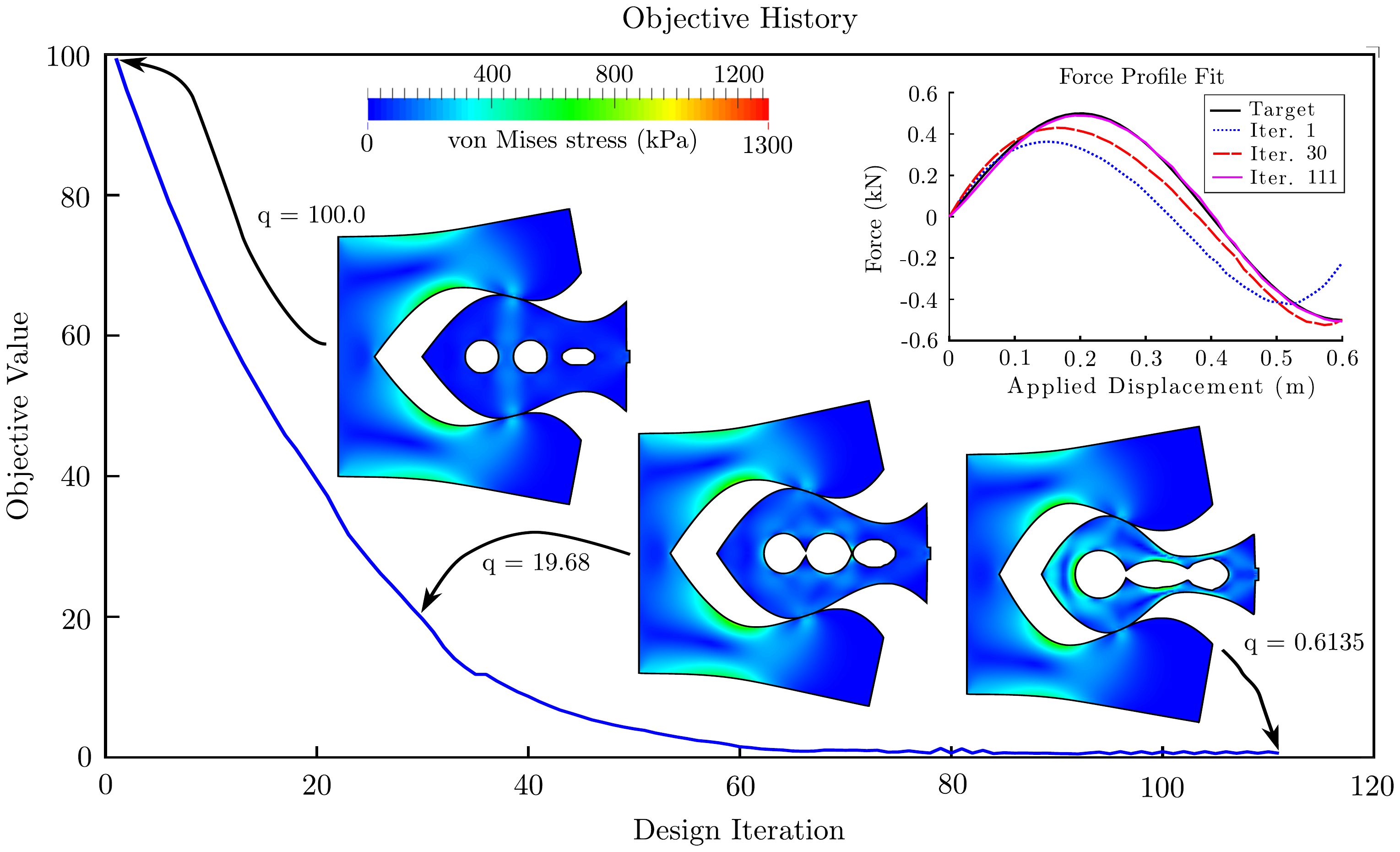}
\caption{Snap fit design objective history with snapshots of specific iterations at applied displacement $U^A_x = 0.213$ m. Inset depicts the force-displacement curve for specific iterations.}
\label{fig:ex2_3p_final}
\end{center}
\end{figure}
%
%\begin{figure}[t]
%\begin{center}
%\includegraphics[width=.5\linewidth]{./figures/ex2_3p_inset.pdf}
%\caption{Snap-fit design force displacement curve, compared to objective.}
%\label{fig:ex2_3p_inset}
%\end{center}
%\end{figure}

The objective history throughout the optimization process is illustrated in Figure \ref{fig:ex2_3p_final}, supported by snapshots of the physical response of the design at $U_x^A = 0.2133$ m for select iterations. The inset of Figure \ref{fig:ex2_3p_final} compares the experienced force-displacement of select design iterations to the desired target profile. Throughout the design evolution, the void material expands until the void inclusions coalesce. The snap-fit general profile narrows, while the peak width increases. The combination of these design attributes affords greater elastic stretch in the snap-fit tab, delaying the maximum force value experienced during incremental loading. This example demonstrates that the proposed optimization method allows finding non-intuitive optimized designs for frictionless contact problems experiencing large deformations.

%===============================================================================

\subsection{Torque Limiter Design Problem}\label{sec:torqueprob}

Torque limiters are common devices used in mechanical equipment to prevent damage from overload. Also known as an overload clutch, these devices limit the applied torque to an assembly by slipping or uncoupling the load. Common methods of limiting applied torque by slipping include frictional plates, magnetic clutches, and ball-detent designs. Taking a simpler approach, the torque limiter design problem presented here consists of two pieces; an outer square shaft containing an inner rod with frictionless contact prescribed at the interface.  The torque limiter design problem is explored for two scenarios: the first example is a two phase design in which geometry control is provided by discretized level set nodal variables; the second example is a three-phase design in which geometry control is provided by geometric primitive variables.

\subsubsection{Two-Phase Example}

\begin{figure}[t]
\begin{minipage}[t]{0.47\textwidth}
\centering
\vspace{0pt}
\includegraphics[width=.8\linewidth]{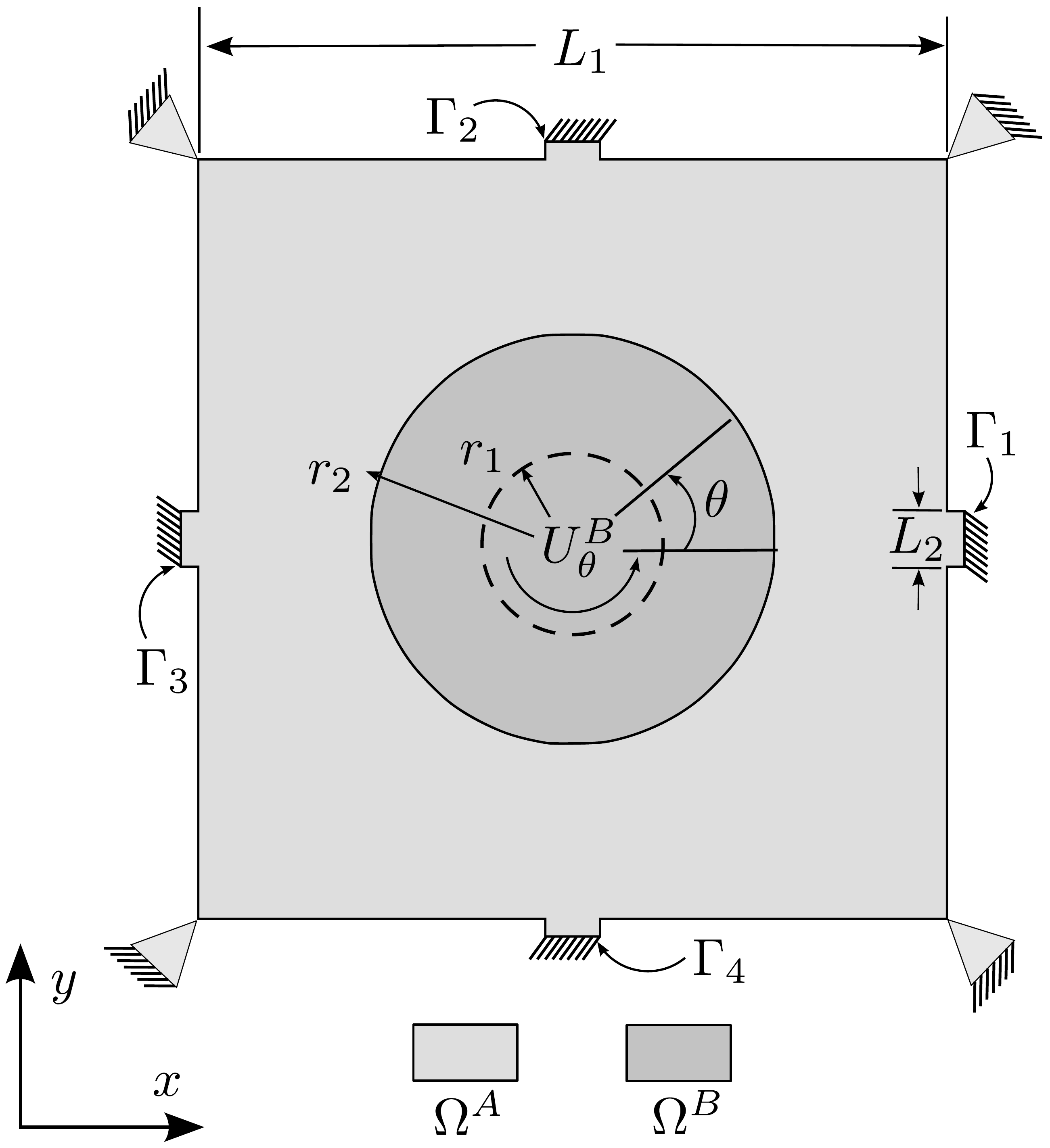}
\caption{Torque limiter initial configuration.}
\label{fig:ex4_setup}
\end{minipage}\hfill
\begin{minipage}[t]{0.51\textwidth}
\centering
\vspace{0pt}
\captionof{table}{Torque limiter model parameters.}
\label{tab:ex4_setup}
\begin{tabular}{ll}
\hline
Description & Parameter \\ \hline
domain length & $L_1 = 1.0$ m \\
fixed support width & $L_2 = 0.0732$ m \\
loading radius & $r_1 = 0.12$ m \\
interface radius & $r_2 = 0.291$ m \\
Young's modulus & $E^A = 10$ Mpa \\
Young's modulus & $E^B = 10$ Mpa \\
Poisson's ratio & $\nu^A = 0.3$ \\
Poisson's ratio & $\nu^B = 0.3$ \\
applied rotation at $t = 1$ & $U_\theta^B = \pi / 2$ rad \\
response weight &  $c_u = 100.0$ \\
penalty weight & $c_p = 0.01$ \\
volume ratio & $c_v = 1$ \\
opt. upper bounds & $s_max = 0.0244$ \\
opt. lower bounds & $s_min = -0.0244$ \\
rel. step size & $\Delta s = 1\times10^{-2}$\\
\hline
\end{tabular}
\end{minipage}
\end{figure}

The torque-limiter design problem consists of an outer shaft, constrained along the outer boundaries, with an internal rod which is rotated as depicted in Figure \ref{fig:ex4_setup}.
The outer square shaft, represented by phase A, is grounded at the corners and at boundaries $\Gamma_{1-4}$. The inner rod of radius $r_2$, represented by phase B, is rotated about the centroid with an applied displacement $U^B_\theta(t)$ within radius $r_1$. The displacement $U^B_\theta(t)$ is applied in $30$ equal load increments. Model parameters specific to this problem are listed in Table \ref{tab:ex4_setup}. The objective of this study is to find the optimal geometry such that the total torque experienced at boundaries $\Gamma_{1-4}$ matches a target torque profile. The objective function is defined as follows:
\begin{equation}\label{eq:ex3_obj}
z = \int \left( \int_{\Gamma_{1-4}}\sigma_{r\theta} \ r \  \mathrm{d}\Gamma  - f_t \right)^2\mathrm{d}t,
\end{equation}
where $\sigma_{r\theta}$ is the shear component of the Cauchy stress in a polar coordinate system, $r$ is the radial position in a polar coordinate system, and $f_t$ is the desired torque profile. For this particular example, the desired torque profile is defined as:
\begin{equation}
f_t = 2.5 \ \sin{2\pi t} \ \mathrm{kN}\cdot \mathrm{m}, \qquad 0\leq t \leq 1.
\end{equation}
The chosen target curve describes a load-displacement curve that follows a sinusoidal wave. Pseudo-time $t$ defines the linear increment of applied rotation during the period of $0\leq t \leq 1$. A sinusoidal target curve provides a gradual transition from experienced torque build-up and decrease, delineated by a peak target value. The entire domain is discretized with 41$\times$41 elements. For geometry control, the nodal level set values are defined as the optimization variables. Numerical experiments showed that while a volume constraint is not necessary for this problem, a perimeter penalty weight of $c_p = 0.01$ is useful to promote a smooth interface profile. Also the smoothing radius of the linear filter \eqref{eq:smooth_rad} is set to $r_f = 1.5\sqrt{2}h$, where $h$ is the element side length for this particular problem.
\begin{figure}[t]
\begin{center}
\includegraphics[width=.95\linewidth]{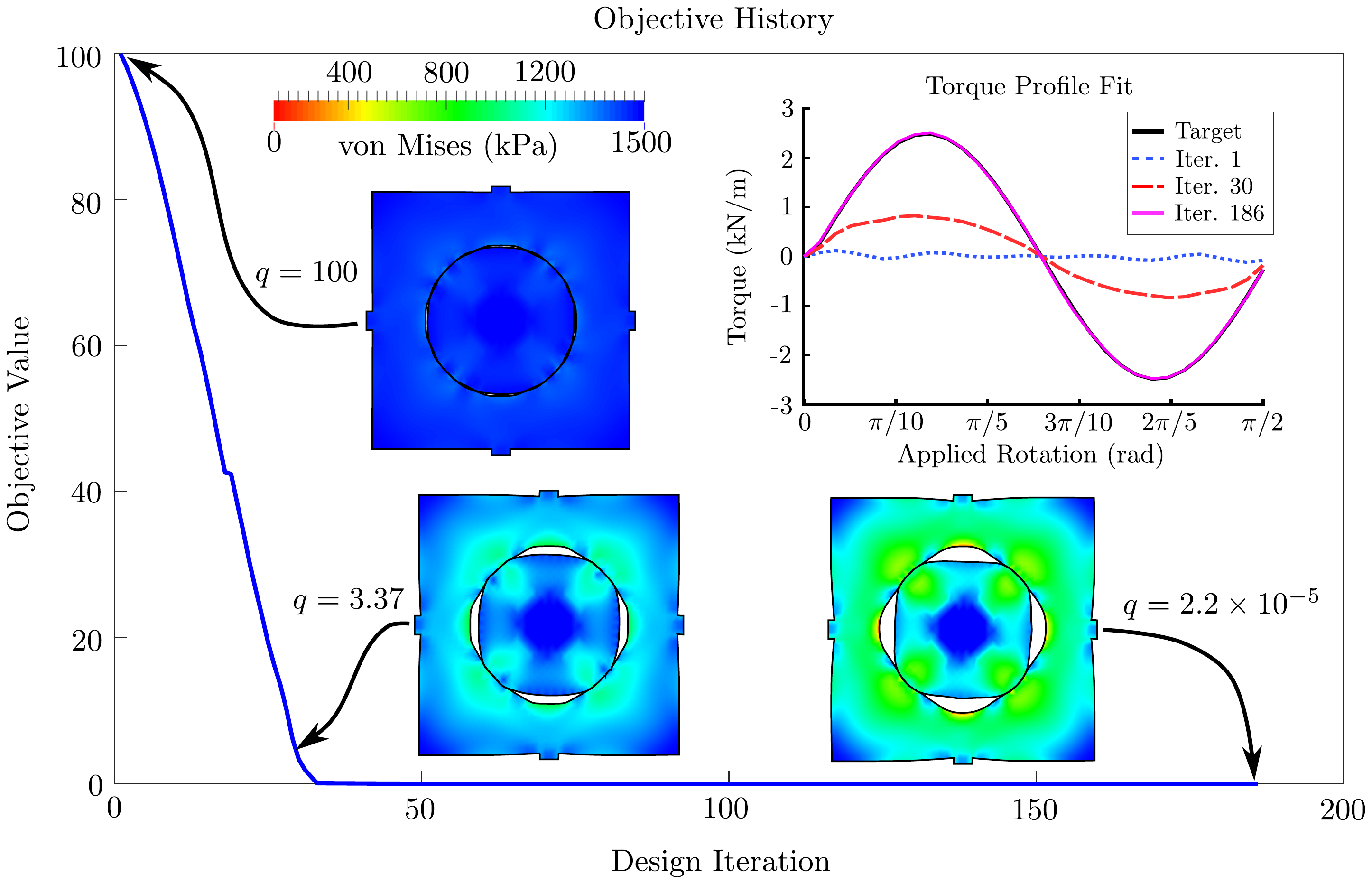}
\caption{Torque limiter design objective history with snapshots of specific iterations. Inset depicts the torque-rotation curve for specific iterations.}
\label{fig:ex4_final}
\end{center}
\end{figure}

Figure \ref{fig:ex4_final} shows the objective history throughout the optimization process, supported by snapshots of the mechanical response at the final time step for select design iterations. The inset of Figure \ref{fig:ex4_final} shows the experienced torque profile, compared to the desired profile for specific design iterations. As expected, the initial circular profile of the inner rod yields no torque as it is rotated. The interface geometry evolves bumps or ridges, providing a torque-rotation profile that closely resembles the desired profile.

\subsubsection{Three-Phase Example}
\begin{figure}[b]
\begin{minipage}[t]{0.47\textwidth}
\centering
\vspace{0pt}
\includegraphics[width=.8\linewidth]{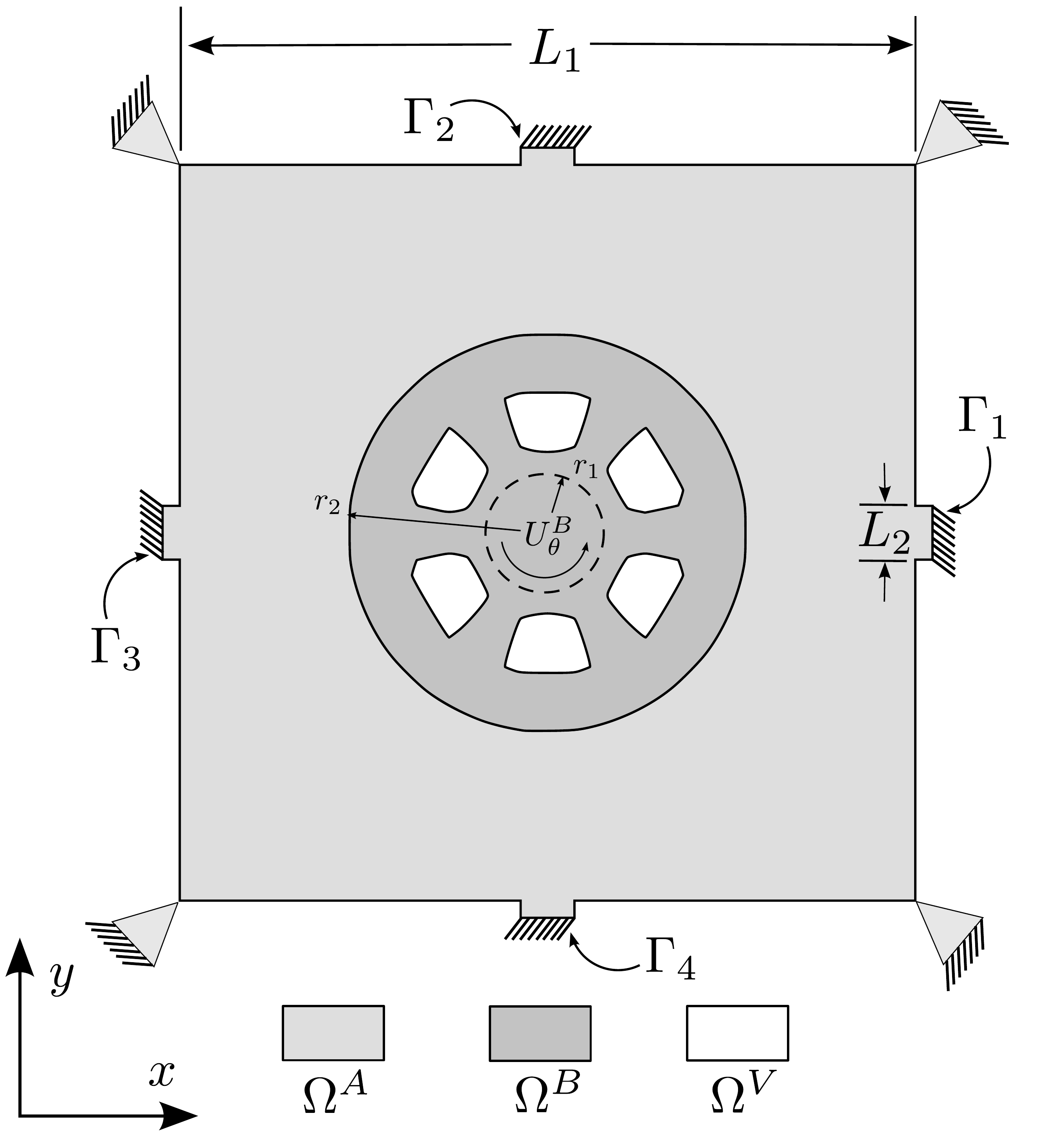}
\caption{Three-phase torque limiter initial configuration.}
\label{fig:ex5_setup}
\end{minipage}\hfill
\begin{minipage}[t]{0.51\textwidth}
\centering
\vspace{0pt}
\captionof{table}{Three-phase torque limiter model parameters.}
\label{tab:ex5_setup}
\begin{tabular}{ll}
\hline
Description & Parameter \\ \hline
domain length & $L_1 = 1.0$ m \\
fixed support width & $L_2 = 0.0732$ m \\
loading radius & $r_1 = 0.06$ m \\
interface radius & $r_2 = 0.291$ m \\
Young's modulus & $E^A = 10$ Mpa \\
Young's modulus & $E^B = 10$ Mpa \\
Poisson's ratio & $\nu^A = 0.3$ \\
Poisson's ratio & $\nu^B = 0.3$ \\
applied rotation at $t = 1$ & $U_\theta^B = \pi / 3$ rad \\
response weight &  $c_u = 99.9$ \\
penalty weight & $c_p = 0.1$ \\
volume ratio & $c_v = 0.16$ \\
rel. step size & $\Delta s = 1\times10^{-2}$\\
\hline
\end{tabular}
\end{minipage}
\end{figure}
The three-phase torque limiter example explores the buildup of strain energy and abrupt release, also known as 'snap through' behavior, by introducing void material within the rod phase, $\Omega^B$. Similar to the two-phase snap-fit design problem of Section \ref{sec:3psfit}, the three-phase torque limiter design problem defines geometry by a set of geometric primitive shapes. The initial configuration is depicted in Figure \ref{fig:ex5_setup}, and model parameters are listed in Table \ref{tab:ex5_setup}.

The initial design geometry closely resembles that of the two-phase torque limiter problem, however, void regions are introduced radially within phase B. Similarly, we wish to determine the best arrangement of material to match the desired torque profile:
\begin{equation}
f_t = 5.20 \ t(1-t)\left(0.5 -  \frac{\tanh(30) + \tanh(60(t - 0.5))}{2\tanh(30)}\right) \ \mathrm{kN}\cdot \mathrm{m}
\qquad 0\leq t \leq 1.
\end{equation}
The torque profile defined for this problem increases to a peak value, then abruptly drops to a minimum value before gradually returning to zero at the final load step. This torque profile encourages a build-up and abrupt release of strain energy. The prescribed displacements are applied in $30$ equal load increments.
\begin{figure}[t]
\begin{center}
\includegraphics[width=.7\linewidth]{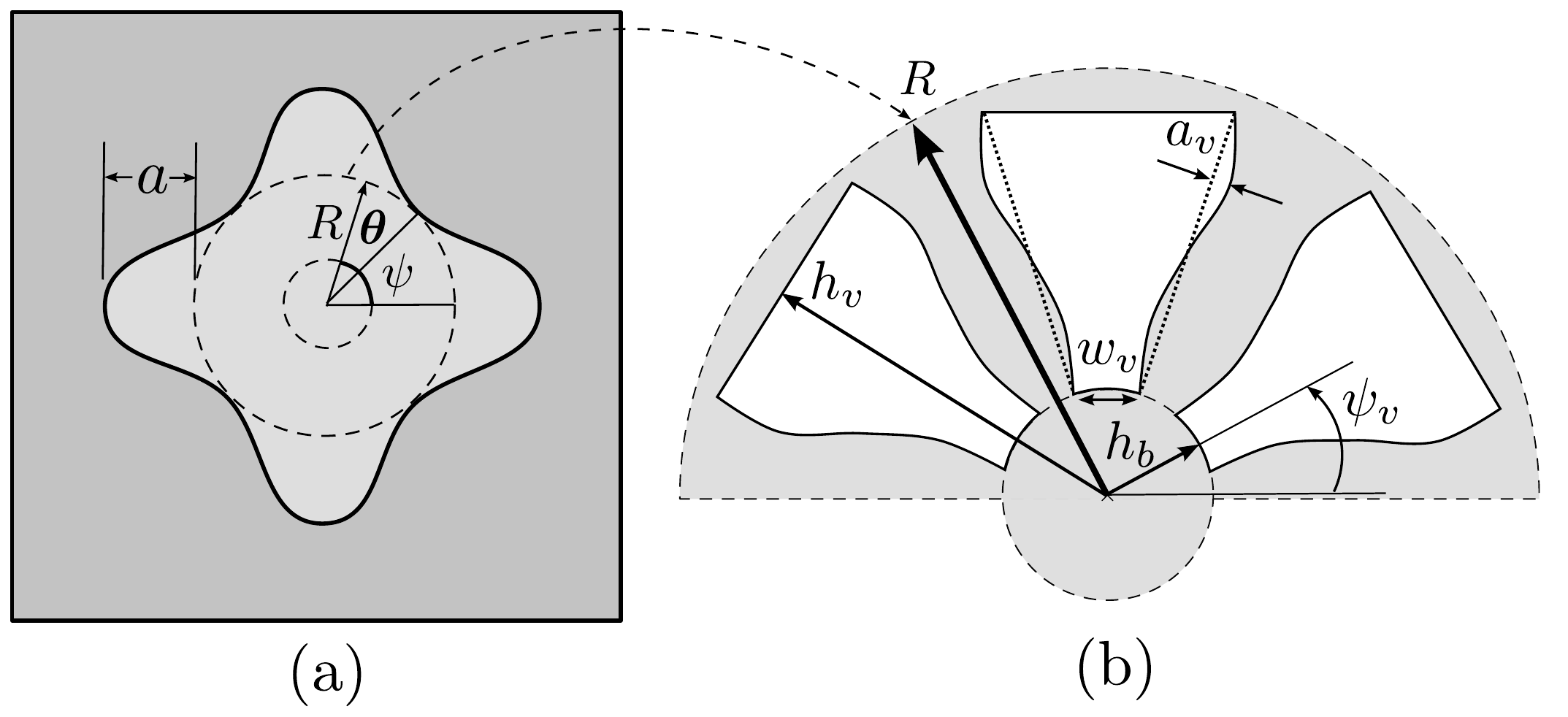}
\caption{Torque limiter geometric primitives for (a)$\ \phi^1$ and (b)$\ \phi^2$.}
\label{fig:ex4_3p_geoprim}
\end{center}
\end{figure}

To afford geometry control, both LSFs define the geometric primitives depicted in Figure \ref{fig:ex4_3p_geoprim}. LSF $\phi^1$ describes a circle with a sinusoidally varying radius:
\begin{equation}
\phi^1 = R - \left| \frac{\left(X-X_c\right)/L^1}{\cos(n \theta)} - a \sin(n \theta + \psi)\right|
\end{equation}
where $R$ is the base circle radius, $a$ is the amplitude of variations, $n$ is the number of ridges along the outer surface, $X_c$ and $Y_c$ are the center of the design domain, $\psi$ is the phase of variations, and the auxiliary coordinate $\theta$ is defined as
\begin{equation}
\theta = \mathrm{atan2}\left(\left(Y-Y_c\right)/L_1, \left(X-x^c\right)/L_1\right),
\end{equation}
where $\mathrm{atan2}$ is the four-quadrant inverse tangent. For LSF $\phi^1$, the number of ridges must be a positive integer value, as non-integer values yield a discontinuous zero level-set contour.

The second LSF, $\phi^2$, is describes a series of void petals as shown in Figure \ref{fig:ex4_3p_geoprim}(b). The LSF describing the petals is taken from \cite{CM1:16}. The $jth$ petal is defined by:
\begin{equation}
\phi_{v,j} = - h_{b} + \left(\left(\frac{2 \tilde{X}_v h_{b}}{\tilde{w}_v} \right)^{10} + \left(\frac{\tilde{Y}_v h_b }{h_{p,j}}\right)^{10} \right)^{1/10}
\end{equation}
The auxiliary coordinates, $\tilde{X}_v$ and $\tilde{Y}_v$, and pedal width, $\tilde{w}_v$, are defined as
\begin{equation}
\tilde{X}_v = \acute{X} - \mathrm{sign}\left(-\acute{X}a_{v,j} \tilde{w}_v \sin\left(\frac{3\pi\left(\acute{Y}- h_b\right)}{2 h_{v,j}}\right) \right),
\end{equation}
\begin{equation}
\tilde{Y}_v = \acute{Y} - h_b,
\end{equation}
\begin{equation}
\tilde{w}_v = w_{v,j} + \pi \frac{\acute{Y}-h_b}{N_p - 1},
\end{equation}
where the number of petals, $N_p$, must be a positive integer. Finally, the rotated coordinate system can be expressed as:
\begin{equation}
\acute{X} = \frac{X - X_c}{L_1}\cos\left(\psi_v\right) - \frac{Y-Y_c}{L_1}\sin\left(\psi_v\right), \qquad \acute{Y} = \frac{X-X_c}{L_1}\sin\left(\psi_v\right) + \frac{Y-Y_c}{L^1}\cos\left(\psi_v\right)
\end{equation}
LSF $\phi^2$ is thus defined as
\begin{equation}
\phi^2 = \min\left(\phi_{v,j}\right)
\end{equation}

For this optimization problem, the number of outer ridges, $n = 4$, center location $X_c = 0.5 m$ and $Y_c = 0.5 m$, and number of void regions, $N_v = 6$, are all held constant. Each void region is distributed radially by increments of $\pi/3$, as shown in Figure \ref{fig:ex5_setup}. The 28 design variable initial values, upper and lower bounds are provided in Table \ref{tab:ex4_3p_vars}. A volume constraint of $c_v = 0.16$ is applied to reduce the overall volume occupied by the internal rod, and a perimeter penalty weight of $c_p = 0.1$ is applied to regularize model geometry. No smoothing filter is used. The design domain is discretized with 41$\times$41 elements.
\begin{table}[t]
\centering
\begin{tabular}{lllll}
\hline
Description & Variable & Initial Value & Upper Bound & Lower Bound \\ \hline
Outer base radius        & $R$           & $0.271$ m & $0.4$ m      & $0.1$ m \\
Outer surface amplitude  & $a$           & $0.0$ m   & $0.05$ m     & $0.0$ m  \\
Phase of outer surface   & $\psi$        & $0.0$ rad & $\pi/2$ rad  & $-\pi/2$ m \\
Petal radial location    & $\psi_{v,j}$  & var.      & $+\Delta\pi/4$ m   & $-\Delta\pi/4$ m \\
Petal base               & $h_b$         & $0.15$ m  & $0.3$ m      & $0.012$ m  \\
Petal height             & $h_{v,j}$     & $0.04$ m  & $0.1$ m      & $0.01$ m  \\
Petal side variation     & $a_{v,j}$     & $0.0$ m   & $0.1$ m      & $-0.1$ m \\
Petal Width              & $w_{v,j}$     & $0.1$ m   & $0.3$ m      & $0.05$ m \\
\hline
\vspace{0.1cm}
\end{tabular}
\caption{Initial value, upper and lower bounds for three-phase torque limiter problem.}\label{tab:ex4_3p_vars}
\end{table}
\begin{figure}[t]
\begin{center}
\includegraphics[width=1\linewidth]{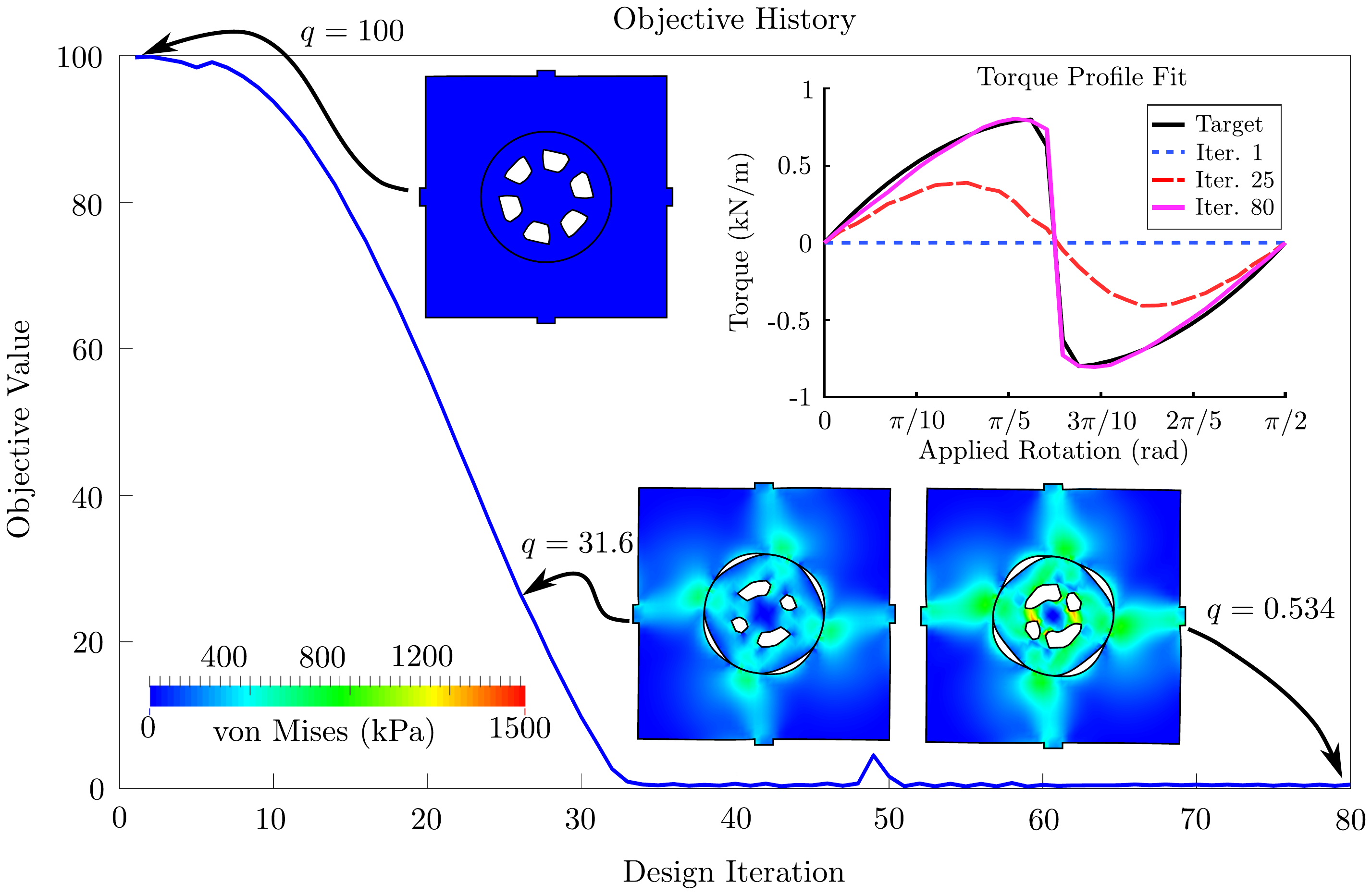}
\caption{Torque limiter design geometry evolution at $U_\theta^B \approx \pi/5$ rad, colored by von Mises stress.}
\label{fig:ex4_3p_hist}
\end{center}
\end{figure}
%
%\begin{figure}[t]
%\begin{center}
%\includegraphics[width=.5\linewidth]{./figures/ex4_3p_inset.pdf}
%\caption{Torque limiter design force displacement curve, compared to target curve.}
%\label{fig:ex4_3p_inset}
%\end{center}
%\end{figure}

Figure \ref{fig:ex4_3p_hist} depicts the objective history of the optimization problem supported by snapshots of select design iterations throughout the optimization process. The inset plot within Figure \ref{fig:ex4_3p_hist} shows the torque-rotation curves for select design iterations compared to the target profile. At early stages of convergence, the contact interface, $\Gamma_c$, evolves to exhibit ridges to increase the experienced torque along the boundaries $\Gamma_{1-4}$. To match the sudden drop in the target torque-rotation curve, the void regions coalesce to reduce the overall material bridging the area of applied displacement to the material in the vicinity of the contact interface. This encourages snap-through behavior. The optimized design matches the desired torque profile well, illustrating the ability of the proposed optimization method to find geometries that feature complex mechanical contact behavior.

%===============================================================================

\section{Conclusions}\label{sec:con}
\vspace{-2pt}
This paper presented a shape and topology optimization framework for two- and three-phase problems with finite strain, large sliding bilateral contact phenomena. Coincident surface location was defined by a coupled parametric representation of the surface geometry. The interface condition was described by a stabilized Lagrange formulation with an active-set strategy to allow surface separation. The material behavior of the mechanical model was described by a hyper-elastic isotropic material and finite strains were assumed for the mechanical model. The XFEM was used to discretize and integrate the mechanical model, and a face-oriented ghost penalization model was used to mitigate the ill-conditioning of the physical response prediction. Dynamic relaxation was employed to provide reliable convergence.

Geometry control was provided by an explicit LSM, where single and multiple LSFs were used to describe two- and three-phase geometries, respectively. The optimization problem was solved with a nonlinear programming method, and a perimeter penalty was used to regularize optimized geometry. Material phase volume constraints were imposed, and design sensitivities were evaluated using an adjoint method. The physical response prediction and subsequent design sensitivities were verified with benchmark examples.

The comparative optimization study between small infinitesimal strain theory and finite strain theory demonstrated that within acceptable load limits for linearized strain theory, both physical response models produced similar results. The two-phase snap-fit and torque limiter optimization studies explored physical response behavior that could not be modeled with infinitesimal strain theory. Furthermore, it was demonstrated that reliable convergence behavior, and subsequently, non-intuitive design solutions are feasible for this particular subset of problems. The three-phase design problems explored in this study demonstrated that void regions within components can contribute to design functionality, specifically by increasing the level of elastic deformation and affording snap through behavior in the snap-fit and torque limiter design problems, respectively.

While this study focused on frictionless and quasi-static contact behavior, the current framework allows for the convenient extension to rate-based interface phenomena such as non-conservative frictional effects. Some of the non-smooth force-displacement behavior may be attributed to the piece-wise linear approximation of interface geometry. Future studies could benefit from $C^1$ continuous discretization. The current optimization method does not allow for the nucleation of a phase within a volume of another phase. Methods for overcoming this limitation, such as the use of topological derivatives, should be investigated. Extension into three dimensional problems would greatly increase the scope of potential applications, as relatively few real world applications can be reduced into two dimensional space. Additional solution techniques such as arc length methods, and optimization methods should be explored to improve the stability of the physical response prediction and optimization convergence behavior. In future studies, this method needs to be extended to the treatment of triple junction intersections and other considerations for a more generalized multi-phase approach.

%===============================================================================

\section*{Acknowledgements}

The authors acknowledge the support of the National Science Foundation under grant CMMI-1235532 and support through Sandia National Laboratories under Contract Agreement 1396470. The second author acknowledges the support of the Air Office of Scientific Research under Contract No. FA9550-13-1-008. This work utilized the Janus supercomputer, which is supported by the National Science Foundation (award number CNS-0821794) and the University of Colorado Boulder. The Janus supercomputer is a joint effort of the University of Colorado Boulder, the University of Colorado Denver and the National Center for Atmospheric Research. Janus is operated by the University of Colorado Boulder. The opinions and conclusions presented in this paper are those of the authors and do not necessarily reflect the views of the sponsoring organization.

\vspace{-6pt}

%===============================================================================

%\bibliographystyle{wileyj}
%\bibliography{JabRef_Database}

\begin{thebibliography}{10}
\providecommand{\url}[1]{\texttt{#1}}
\providecommand{\urlprefix}{URL }
\expandafter\ifx\csname urlstyle\endcsname\relax
  \providecommand{\doi}[1]{doi:\discretionary{}{}{}#1}\else
  \providecommand{\doi}{doi:\discretionary{}{}{}\begingroup
  \urlstyle{rm}\Url}\fi

\bibitem{HKP:99}
Hilding D, Klarbring A, Petersson J. Optimization of structures in unilateral
  contact. \emph{Appl. Mech. Rev.}  1999; \textbf{52}(4):139--160.

\bibitem{HTK:01}
Hilding D, Torstenfelt B, Klarbring A. A computational methodology for shape
  optimization of structures in frictionless contact. \emph{Comput. Methods
  Appl. Mech. Engrg.}  2001; \textbf{190}:4043--4060.

\bibitem{HLDS:00}
Herskovits J, Leontiev A, Dias G, Santos G. Contact shape optimization: a
  bilevel programming approach. \emph{Struct Multidisc Optim}  2000;
  \textbf{20}:214--221.

\bibitem{SMR:00}
Schleupen A, Maute K, Ramm E. Adaptive fe-procedures in shape optimization.
  \emph{Structural and Multidisciplinary Optimization}  2000;
  \textbf{19}:282--302.

\bibitem{Bendsoe:89}
Bends{\o}e M. Optimal shape design as a material distribution problem.
  \emph{Structural and Multidisciplinary Optimization}  1989;
  \textbf{1}(4):193--202, \doi{10.1007/BF01650949}.

\bibitem{RZB:92}
Rozvany G, Zhou M, Birker T. Generalized shape optimization without
  homogenization. \emph{Structural and Multidisciplinary Optimization}  1992;
  \textbf{4}(3):250--252.

\bibitem{BS:03}
Bends{\o}e MP, Sigmund O. \emph{Topology Optimization: Theory, Methods and
  Applications}. Springer, 2003.

\bibitem{SM:13}
Sigmund O, Maute K. Topology optimization approaches: {A} comparative review.
  \emph{Structural and Multidisciplinary Optimization}  2013;
  \textbf{48}(6):1031--1055.

\bibitem{DG:14}
Deaton JD, Grandhi RV. A survey of structural and multidisciplinary continuum
  topology optimization: post 2000. \emph{Structural and Multidisciplinary
  Optimization}  2014; \textbf{49}(1):1--38, \doi{10.1007/s00158-013-0956-z}.

\bibitem{KC:01a}
Kikuchi N, Chen BC. Topology optimization with design-dependent loads.
  \emph{Finite Elements in Analysis and Design}  2001; \textbf{37}(1):57--70.

\bibitem{SC:07}
Sigmund O, Clausen P. Topology optimization using a mixed formulation: An
  alternative way to solve pressure load problems. \emph{Computer Methods in
  Applied Mechanics and Engineering}  2007; \textbf{196}(13--16):1874--1889,
  \doi{http://dx.doi.org/10.1016/j.cma.2006.09.021}.

\bibitem{Yoon:10}
Yoon GH. Topology optimization for stationary fluid--structure interaction
  problems using a new monolithic formulation. \emph{International Journal for
  Numerical Methods in Engineering}  2010; \textbf{82}(5):591--616.

\bibitem{HO:00}
Hammer V, Olhoff N. Topology optimization of continuum structures subjected to
  pressure loading. \emph{Struct. Multidisc. Optim.}  2000; \textbf{19}:85--92.

\bibitem{ARS:12}
Andrade-Campos A, Ramos A, Sim$\tilde{\mathrm{o}}$es J. A model of bone
  adaptation as a topology optimization process with contact. \emph{Journal of
  Biomedical Science and Engineering}  2012; \textbf{5}:229--244.

\bibitem{S:09}
Str\"{o}mberg N. The influence of sliding friction in optimal topologies.
  \emph{Recent Advances in Contact Mechanics}, vol.~56, Stravroulakis GE (ed.),
  5th Contact Mechanics International Symposium (CMIS2009), Springer, 2009;
  327--336.

\bibitem{LLK:15}
Luo Y, Li M, Kang Z. Topology optimization of hyperelastic structures with
  frictionless contact supports. \emph{International Journal of Solids and
  Structures}  2015; \textbf{81}:373--382.

\bibitem{DML+:13}
van Dijk N, Maute K, Langelaar M, Keulen F. Level-set methods for structural
  topology optimization: a review. \emph{Structural and Multidisciplinary
  Optimization}  2013; :1--36\doi{10.1007/s00158-013-0912-y}.

\bibitem{M:09}
My\'sli\'nski A. Level set method for shape and topology optimization of
  contact problems. \emph{System Modeling and Optimization}, \emph{IFIP
  Advances in Information and Communication Technology}, vol. 312, Korytowski
  A, Malanowski K, Mitkowski W, Szymkat M (eds.), Springer Berlin Heidelberg,
  2009; 397--410.

\bibitem{LM:15}
Lawry M, Maute K. Level set topology optimization of problems with sliding
  contact interfaces. \emph{Structural and Multidisciplinary Optimization}
  2015; \textbf{52}(6):1107--1119.

\bibitem{LLK:16}
Liu P, Luo Y, Kang Z. Multi-material topology optimization considering
  interface behavior via {XFEM} and level set method. \emph{Computer Methods in
  Applied Mechanics and Engineering}  2016; \textbf{308}:113--133.

\bibitem{FB:06}
Fries TP, Belytschko T. The intrinsic xfem: A method for arbitrary
  discontinuities without additional unkowns. \emph{International Journal for
  Numerical Methods in Engineering}  2006; \textbf{68}:1358--1385,
  \doi{10.1002/nme.1761}.

\bibitem{K:15}
Khoei AR. \emph{Extended finite element method: Theory and applications}.
  Wiley, 2015.

\bibitem{LB:08}
Liu F, Borja R. A contact algorithm for frictional crack propagation with the
  extended finite element method. \emph{International Journal For Numerical
  Methods In Engineering}  2008; \textbf{76}:1489--1512.

\bibitem{MWL:12}
Mueller-Hoeppe DS, Wriggers P, Loehnert S. Crack face contact for
  hexahedral-based {XFEM} formulation. \emph{Computational Mechanics}  2012;
  \textbf{49}:725--734.

\bibitem{GMM:07}
Geniaut S, Massin P, Mo{\"e}s N. A stable {3D} contact formulation using
  {X-FEM}. \emph{REMN}  2007; \textbf{16}:259--275.

\bibitem{AHL:12}
Amdouni S, Hild P, Lleras V, Moakher M, Renard Y. A stabilised {L}agrange
  multiplier method for the enriched finite-element approximation of contact
  problems of cracked elastic bodies. \emph{ESAIM: Mathematical Modelling and
  Numerical Analysis}  2012; \textbf{46}:813--839.

\bibitem{LB:10}
Liu F, Borja RI. Stabilized low-order finite elements for frictional contact
  with the extended finite element method. \emph{Computer Methods in Applied
  Mechancis and Engineering}  2010; \textbf{199}:2456--2471.

\bibitem{GTT+:10}
Giner E, Tur M, Tarancon JE, Fuenmayor FJ. Crack face contact in {X-FEM} using
  a segment-to-segment approach. \emph{International Journal for Numercial
  Methods in Engineering}  2010; \textbf{82}:1424--1449.

\bibitem{SGM+:13}
Siavelis M, Guiton M, Massin P, M\"{o}es N. Large sliding contact along
  branched discontinuities with x-fem. \emph{Comput Mech}  2013;
  \textbf{52}:201--219.

\bibitem{NGM+:09}
Nistor I, Guiton MLE, Massin P, Mo\"{e}s N, G\'{e}niaut S. An {X-FEM} approach
  for large sliding contact along discontinuities. \emph{International Journal
  For Numerical Methods In Engineering}  2009; \textbf{78}:1407--1435.

\bibitem{BP:12}
Biotteau E, Ponthot JP. Modeling frictional contact conditions with the penalty
  method in the extended finite element framework. \emph{European Congress on
  Computational Methods in Applied Sciences and Engineering}, 2012.

\bibitem{TM:11}
{Taheri Mousavi} SMJ, {Taheri Mousavi} SM. Modeling large sliding frictional
  contact along non-smooth discontinuities in {X-FEM}. \emph{International
  Journal of Modeling and Optimization}  2011; \textbf{1}:169--173.

\bibitem{WWG:03}
Wang MY, Wang X, Guo D. A level set method for structural topology
  optimization. \emph{Computer Methods in Applied Mechanics and Engineering}
  2003; \textbf{192}(1-2):227--246, \doi{10.1016/S0045-7825(02)00559-5}.

\bibitem{AJT:04}
Allaire G, Jouve F, Toader AM. Structural optimization using sensitivity
  analysis and a level-set method. \emph{Journal of Computational Physics}
  2004; \textbf{194}(1):363--393, \doi{10.1016/j.jcp.2003.09.032}.

\bibitem{KM:12}
Kreissl S, Maute K. Levelset based fluid topology optimization using the
  extended finite element method. \emph{Structural and Multidisciplinary
  Optimization}  2012; :1--16.

\bibitem{VC:02}
Vese L, Chan T. A multiphase levelset framework for image segmentation using
  the mumford and shah model. \emph{International Journal of Computer Vision}
  2002; \textbf{50}(3):271--293.

\bibitem{DAP+:12}
Dombre E, Allaire G, Pantz O, Schmitt D. Shape optimization of a sodium fast
  reactor core. \emph{ESAIM proceedings, EDP Sciences}  2012;
  \textbf{38}:319--334.

\bibitem{WW:05}
Wang M, Wang X. A level-set based variational method for design and
  optimization of heterogeneous objects. \emph{Computer-Aided Design}  2005;
  \textbf{37}(3):321--337.

\bibitem{YX:04a}
Yulin M, Xiaoming W. A level set method for structural topology optimization
  and its applications. \emph{Advances in Engineering Software}  2004;
  \textbf{35}(7):415--441.

\bibitem{Wriggers:02}
Wriggers P. \emph{Computational Contact Mechanics}. Wiley, 2002.

\bibitem{MM:13}
Makhija D, Maute K. Numerical instabilities in level set topology optimization
  with the extended finite element method. \emph{Structural and
  Multidisciplinary Optimization}  2014; \textbf{49}(2):185--197.

\bibitem{TAY:03}
Terada K, Asai M, Yamagishi M. Finite cover method for linear and non-linear
  analyses of heterogeneous solids. \emph{International Journal for Numerical
  Methods in Engineering}  2003; \textbf{58}(9):1321--1346,
  \doi{10.1002/nme.820}.

\bibitem{TYHT:11}
Tran AB, Yvonnet J, He QC, Toulemonde C, Sanahuja J. A multiple level set
  approach to prevent numerical artefacts in complex microstructures with
  nearby inclusions within {XFEM}. \emph{International Journal for Numerical
  Methods in Engineering}  2011; \textbf{85}(11):1436--1459.

\bibitem{AHD:14}
Annavarapu C, Hautefeuille M, Dolbow JE. A {N}itsche stabilized finite element
  method for frictional sliding on embedded interfaces. {P}art {I}: Single
  interface. \emph{Computer Methods in Applied Mechanics and Engineering}
  2014; \textbf{268}:417--436.

\bibitem{BH:12}
Burman E, Hansbo P. Fictitious domain finite element methods using cut
  elements: {II}.~ {A} stabilized {N}itsche method. \emph{Applied Numerical
  Mathematics}  2012; \textbf{62}(4):328--341.

\bibitem{BFH:06}
Burman E, Fern\'{a}ndez M, Hansbo P. Continuous interior penalty finite element
  method for oseen's equations. \emph{SIAM journal on numerical analysis}
  2006; \textbf{44}(3):1248--1274.

\bibitem{SW:14}
Schott B, Wall W. A new face-oriented stabilized {XFEM} approach for {2D} and
  {3D} incompressible navier--stokes equations. \emph{Computer Methods in
  Applied Mechanics and Engineering}  2014; \textbf{276}:233--265.

\bibitem{SRG+:14}
Schott B, Rasthofer U, Gravemeier V, Wall W. A face-oriented stabilized
  {N}itsche-type extended variational multiscale method for incompressible
  two-phase flow. \emph{International Journal for Numerical Methods in
  Engineering}  2014; .

\bibitem{M:78}
Mor\'{e} J. The {L}evenberg-{M}arquardt algorithm: implementation and theory.
  \emph{Numerical Analysis}  1978; :105--116.

\bibitem{K:09}
Kawamoto A. Stabilization of geometrically nonlinear topology optimization by
  the {L}evenberg-{M}arquardt method. \emph{Struct. Multidisc. Optim.}  2009;
  \textbf{37}:429--433.

\bibitem{Svanberg:02}
Svanberg K. A class of globally convergent optimization methods based on
  conservative convex separable approximations. \emph{SIAM J. on Optimization}
  2002; \textbf{12}(2):555--573, \doi{10.1137/s1052623499362822}.

\bibitem{AK:96}
Antunez H, Kleiber M. Sensitivity analysis of metal forming process involving
  frictional contact in steady state. \emph{Materials Processing Technology}
  1996; \textbf{60}(1–4):485–491.

\bibitem{KPC:01}
Kim N, Park YH, Choi K. Optimization of a hyperelastic structure with multibody
  contact using continuum-based shape design sensitivity analysis.
  \emph{Struct. Optim.}  2001; \textbf{21}(3):196–208.

\bibitem{KYC:02}
Kim NH, Yi K, Choi KK. A material derivative approach in design sensitivity
  analysis of three-dimensional contact problems. \emph{International Journal
  of Solids and Structures}  2002; \textbf{39}(8):2087 -- 2108.

\bibitem{CM1:16}
Coffin P, Maute K. Level set topology optimization of cooling and heating
  devices using a simplified convection model. \emph{Structural and
  Multidisciplinary Optimization}  2016; \textbf{53}(5):985 -- 1003.

\end{thebibliography}

%===============================================================================

\end{document}